\documentclass[oneside,leqno,11pt]{article}
\usepackage{amssymb,amsmath,amscd,latexsym}

\def\ga{\mathfrak{a}}

\def\gg{\mathfrak{g}}
\def\gh{\mathfrak{h}}

\def\gk{\mathfrak{k}}

\def\gm{\mathfrak{m}}
\def\gn{\mathfrak{n}}
\def\go{\mathfrak{o}}
\def\gp{\mathfrak{p}}

\def\gs{\mathfrak{s}}
\def\gt{\mathfrak{t}}
\def\gu{\mathfrak{u}}
\def\gv{\mathfrak{v}}
\def\gw{\mathfrak{w}}

\def\gz{\mathfrak{z}}

\def\ul{\mathbf{l}}
\def\uell{\mathbb{\ell}}
\def\um{\mathbf{m}}

\def\uw{\mathbf{w}}

\def\uz{\mathbf{z}}

\def\ggg{> \hskip -5 pt >}

\def\Aut{{\rm Aut}}

\def\Ad{{\rm Ad}}

\def\rank{{\rm rank\,}}

\def\Ind{{\rm Ind\,}}

\def\Im{{\rm Im\,}}
\def\Re{{\rm Re\,}}
\def\Pf{{\rm Pf}}

\def\Skew{{\rm Skew\,}}  

\def\C{\mathbb{C}}

\def\F{\mathbb{F}}

\def\H{\mathbb{H}}

\def\O{\mathbb{O}}

\def\R{\mathbb{R}}

\def\cA{\mathcal{A}}

\def\cC{\mathcal{C}}

\def\cE{\mathcal{E}}
\def\cF{\mathcal{F}}

\def\cH{\mathcal{H}}

\def\cM{\mathcal{M}}

\def\cO{\mathcal{O}}

\def\cU{\mathcal{U}}

\newtheorem{theorem}[equation]{Theorem}

\newtheorem{lemma}[equation]{Lemma}

\newtheorem{corollary}[equation]{Corollary}

\newtheorem{proposition}[equation]{Proposition}

\newtheorem{definition}[equation]{Definition}

\title{Infinite Dimensional Multiplicity Free Spaces III:\\
Matrix Coefficients and Regular Functions}

\author{Joseph A. Wolf\footnote{Research partially supported by the 
National Science Foundation
\endgraf
\endgraf
}}

\date{8 September 2009}

\begin{document}


\maketitle

\begin{abstract}
In earlier papers we studied 
direct limits $(G,K) = \varinjlim\, (G_n,K_n)$ of two types of Gelfand 
pairs.  The first type was that in which the $G_n/K_n$ are compact riemannian
symmetric spaces.  The second type was that in which $G_n = N_n\rtimes K_n$ 
with $N_n$ nilpotent,
in other words pairs $(G_n,K_n)$ for which $G_n/K_n$ is a commutative
nilmanifold.  In each we worked out a method inspired by the
Frobenius--Schur Orthogonality Relations to define isometric injections
$\zeta_{m,n}: L^2(G_n/K_n) \hookrightarrow L^2(G_m/K_m)$ for $m \geqq n$ 
and prove that the left regular representation of $G$ on the Hilbert space 
direct limit $L^2(G/K) := \varinjlim L^2(G_n/K_n)$ is multiplicity--free.
This left open questions concerning the nature of the elements of
$L^2(G/K)$.  Here we define spaces $\cA(G_n/K_n)$ of regular functions on
$G_n/K_n$ and injections $\nu_{m,n} : \cA(G_n/K_n) \to \cA(G_m/K_m)$
for $m \geqq n$ related to restriction by $\nu_{m,n}(f)|_{G_n/K_n} = f$.
Thus the direct limit $\cA(G/K):= \varinjlim \{\cA(G_n/K_n), \nu_{m,n}\}$
sits as a particular $G$--submodule of the much larger inverse limit
$\varprojlim \{\cA(G_n/K_n), \text{restriction}\}$.
Further, we define a pre Hilbert space structure
on $\cA(G/K)$ derived from that of $L^2(G/K)$.  This allows an interpretation
of $L^2(G/K)$ as the Hilbert space completion of the concretely defined 
function space $\cA(G/K)$, and also defines
a $G$--invariant inner product on $\cA(G/K)$ for which the left regular 
representation of $G$ is multiplicity--free.  
\end{abstract}

\section{Introduction} \label{sec1}
\setcounter{equation}{0}
Gelfand pairs $(G,K)$, and the corresponding commutative homogeneous
spaces $G/K$ generalize both the concept of riemannian symmetric space 
and that of a familiar family of riemannian nilmanifolds.
Let $G$ be a locally compact topological group and $K$ a compact subgroup.
Then $(G,K)$ is a Gelfand pair, i.e. $L^1(K\backslash G/K)$ is
commutative under convolution, if and only if
the (left regular) representation of $G$ on $L^2(G/K)$ is multiplicity free.

In two earlier notes (\cite{W4} and \cite{W5}) we looked at certain
cases where $G$ and $K$ are not locally compact, in fact are infinite 
dimensional of the form $\varinjlim (G_n,K_n)$ where the $(G_n,K_n)$ are
finite dimensional Gelfand pairs in the usual sense.  We showed in those 
cases that the multiplicity--free condition  is satisfied when
$L^2(G/K)$ is interpreted as a certain Hilbert space direct limit of the
usual $L^2(G_n/K_n)$.  That direct limit made use of renormalizations
inherent in the Frobenius--Schur othogonality relations for matrix coefficients
of irreducible unitary representations, so the meaning of the elements of
$L^2(G/K)$ was not immediate.  
\medskip

Here we define a ring of regular functions
$\cA(G/K)$ whose nature is transparent.  It is
a subalgebra of $\varprojlim \cA(G_n/K_n)$ where $\cA(G_n/K_n)$ is the
ring of regular functions in the usual sense and $\varprojlim \cA(G_n/K_n)$
is defined by restriction of functions.  This is accompanied by a
$G$--equivariant injective map $\cA(G/K) \to L^2(G/K)$ with dense 
image.  That leads to the main results of this note, both for the direct limits
of compact symmetric spaces of \cite{W4} in Section \ref{sec3} 
and for the direct limits of commutative nilmanifolds of \cite{W5} in
Sections \ref{sec4} and \ref{sec5}.  In both cases the results are that

\begin{equation}\label{generalresults}
\begin{aligned}
&\text{$\cA(G/K)$ injects to a dense subspace of $L^2(G/K)$, so $L^2(G/K)$
defines a $G$--invariant} \\
&\text{inner product on $\cA(G/K)$, the regular 
representation of $G$ on $\cA(G/K)$ is unitarized, } \\
&\text{and $L^2(G/K)$ can be interpreted as the Hilbert space 
completion of $\cA(G/K)$.}
\end{aligned}
\end{equation}

\section{The Ring of Regular Functions}\label{sec2}
\setcounter{equation}{0}

In this section we describe the general setup needed for constructing 
the ring of regular functions on direct systems of Lie groups and 
commutative spaces, and compare them with $L^2$ direct limits. 
We will specialize to parabolic direct systems of compact Lie groups 
compact riemannian symmetric spaces in Section \ref{sec3}.  Then we further 
specialize the results of 
this section to direct systems of certain classes of nilpotent Lie groups 
and commutative nilmanifolds in Sections \ref{sec4} and \ref{sec5}.
\medskip

\subsection{Direct Limit Representations.}\label{sec2a}
We consider direct limit groups $G = \varinjlim G_n$ and direct limit
unitary representations $\pi = \varinjlim \pi_n$ of them.  This means that
$\pi_n$ is a unitary representation of $G_n$ on a Hilbert space $\cH_{\pi_n}$, 
that the $\cH_{\pi_n}$ form a direct system whose maps are unitary, and that 
$\pi$ is the representation of $G$ on $\cH_\pi =\varinjlim \cH_{\pi_n}$ 
given by 
\begin{equation}
\pi(g)v = \pi_n(g_n)v_n \text{ for } n\ggg 0 \text{ so that }
V_n \hookrightarrow V \text{ and } G_n \hookrightarrow
G \text{ send } v_n \to v \text{ and } g_n \to g.
\end{equation}
This formal definition amounts to saying that $\pi$ is a well defined unitary
representation of $G$ on $\cH_\pi$.
\medskip

If $u,v \in \cH_{\pi_n}$ we have the {\sl matrix coefficient}
\begin{equation}\label{def-coef}
f_{u,v,n}: G_n \to \C \text{ defined by }
	f_{u,v,n}(g) = \langle u , \pi_n(g)v\rangle_{\cH_{\pi_n}}.
\end{equation}
Since $\cH_{\pi_n} \to \cH_{\pi_m}$ is a $G_n$--equivariant unitary
injection for $m \geqq n$, we may view $\cH_{\pi_n}$ as a subspace of 
$\cH_{\pi_m}$.  This done, we have
\begin{equation}\label{restr-coef}
\langle u , \pi_n(g)v\rangle_{\cH_{\pi_n}} =
\langle u , \pi_m(g)v\rangle_{\cH_{\pi_m}} \text{ for } u, v \in \cH_{\pi_n}
\text{ and } g \in G_n\,,
\end{equation}
in other words $f_{u,v,n} = f_{u,v,m}|_{G_n}$ for $u,v \in \cH_{\pi_n}$.
Now these coefficients form a direct system.  Formally, we have
\begin{equation}\label{lim-coef}
\begin{aligned}
&\cA(\pi_n) := \{\text{finite linear combinations of the } f_{u,v,n}
\text{ where } u, v \in \cH_{\pi_n}\} \\
&\cA(\pi_n) \hookrightarrow \cA(\pi_m) \text{ using (\ref{restr-coef}), and }
\cA(\pi) := \varinjlim \cA(\pi_n).
\end{aligned}
\end{equation}
In words, $\cA(\pi_n)$ is the space of {\sl regular functions} on $G_n$
defined by $\pi_n$ and $\cA(\pi)$ is the space of {\sl regular functions}
on $G$ defined by $\pi$.  They correspond to the idea of trigonometric
polynomials inside $L^2$ of the circle, or more generally to the idea of
Harish--Chandra module.
\medskip

Note that $\cA(\pi)$ is contained in the 
projective limit $\varprojlim \cA(\pi_n)$ (defined by restriction of
functions), but it is much smaller.

\subsection{Square Integrable Direct Limits}\label{sec2b}
From now on we assume that the groups $G_n$ are
separable, unimodular, locally compact and of Type I. 
(We will be dealing with commutative spaces $G_n/K_n$,
and the commutativity implies unimodularity for $G_n$).  Then one
has the classical decomposition
\begin{equation}\label{planch--exp-finite}
L^2(G_n) = \int_{\widehat{G_n}} (\cH_{\pi_n}\widehat{\otimes}\cH_{\pi_n}^*)
        d(\pi_n)
\end{equation}
where $d(\pi_n)$ is Plancherel measure on the unitary dual $\widehat{G_n}$ and
the projective tensor product $\cH_{\pi_n}\widehat{\otimes}\cH_{\pi_n}^*$
contains $\cA(\pi_n)$ as a dense subspace.  This expresses a function
$f \in L^1(G_n) \cap L^2(G_n)$ in the form
$f(g) = \int_{\widehat{G_n}} \pi_n(f) d(\pi_n)$ where $\pi_n(f) =
\int_{G_n} f(g)\pi_n(g) dg \in \cH_{\pi_n}\widehat{\otimes}\cH_{\pi_n}^*$
is a Hilbert--Schmidt operator on $\cH_{\pi_n}$ and where one has the
Plancherel formula
$||f||_2^2 = \int_{\widehat{G_n}} ||\pi_n(f)||_{HS}^2 d(\pi_n)$.
\medskip

Fix a closed subgroup $Z_n$ in $G_n$ that is co--compact in the center.
It can be $\{1\}$ if $G_n$ has compact center.  If $\chi \in \widehat{Z_n}$
let $\widehat{G_{n,\chi}} = \{\pi_n \in \widehat{G_n} \mid \pi_n(xz)
= \chi(z)^{-1}\pi_n(x)\}$ and consider the Hilbert space 
$L^2(G_{n,\chi})$ of $L^2$ sections of the 
$G_n$--bundle $(G_n \times_\chi \C) \to G_n/Z_n$.  Then
\begin{equation}\label{cent-decomp}
L^2(G_{n,\chi}) = \int_{\widehat{G_{n,\chi}}} (\cH_{\pi_n}
	\widehat{\otimes}\cH_{\pi_n}^*) d_{n,\chi}(\pi_n)
\end{equation}
where $d_{n,\chi}(\pi_n)$ is Plancherel measure on $\widehat{G_{n,\chi}}$.
The decomposition (\ref{planch--exp-finite}) then decomposes further into
\begin{equation}\label{planch--exp-finite-rel}
L^2(G_n) = \int_{\widehat{Z_n}} L^2(G_{n,\chi}) d_n(\chi)
\end{equation}

An irreducible unitary 
representation $\pi_n \in \widehat{G_n}$ is called {\it square integrable}
if its coefficients $f_{u,v,n}$ satisfy $|f_{u,v,n}| \in L^2(G_n/Z_n)$.
This makes sense because
$Z_n$ is co--compact in the center of $G_n$.  Then one has the 
Godement--Frobenius--Schur orthogonality relations.  In particular there 
is a number $\deg \pi_n > 0$ (called the {\it formal degree}) such that
\begin{equation}\label{orthog}
||f_{u,v,n}||^2_{L^2(G_n/Z_n)} = \frac{1}{\deg \pi_n}
||u||^2_{\cH_{\pi_n}} ||v||^2_{\cH_{\pi_n}}\, .
\end{equation}

Now suppose that we have the $\widehat{G_n}$ arranged so that
we have direct systems $\{\pi_n\}$ of unitary representations
and, for each $n$, every $\pi_n \in \widehat{G_n}$ belongs to 
exactly one of those direct
systems.  Write $\beta_{m,n,\pi} : \cH_{\pi_n} \to \cH_{\pi_m}$ for the
direct system maps on the representation spaces.
For each such direct system we have the direct limit
representation $\pi = \varinjlim \pi_n$ of $G = \varinjlim G_n$ 
and the direct limit representation space 
$\cH_\pi = \varinjlim \{\cH_{\pi_n},\beta_{m,n,\pi}\}$.  
\medskip

We now make two strong assumptions on the system $\{G_n\}$.
First, suppose that we can (and do) choose the co--compact closed central
subgroups $Z_n \subset G_n$ so that 
\begin{equation}\label{central-gen}
G_n \hookrightarrow G_{n+1} \text{ maps } Z_n \cong Z_{n+1}\, .
\end{equation}
Then we can write $Z$ for each of the groups $Z_n$, and each of the direct
systems $\{\pi_n\}$ has a common central character $\chi \in \widehat{Z}$.
We now make the further assumption on the groups $G_n$ that
\begin{equation}\label{req-sq-int}
\text{Plancherel--almost--all } \pi_n \in \widehat{G_n}
\text{ are square integrable: their } |f_{u,v,n}| \in
L^2(G_n/Z).
\end{equation}
As described in \cite{W4} and \cite{W5} we use (\ref {orthog}) 
and (\ref{req-sq-int}) to 
scale the inclusions $\cH_{\pi_n}\widehat{\otimes}\cH_{\pi_n}^* 
\hookrightarrow \cH_{\pi_m}\widehat{\otimes}\cH_{\pi_m}^*$ 
for $m \geqq n$ by means 
of formal degrees and obtain $(G_n\times G_n)$--equivariant
isometric inclusions
\begin{equation}\label{incl-gen}
\zeta_{m,n,\pi}: \cH_{\pi_n}\widehat{\otimes}\cH_{\pi_n}^* \to
\cH_{\pi_m}\widehat{\otimes}\cH_{\pi_m}^* 
\text{ defined by }
f_{u,v,n} \mapsto \left ( \tfrac{\deg \pi_m}{\deg \pi_n}\right )^{1/2}f_{u,v,m}.
\end{equation}
In other words $\zeta_{m,n,\pi} = \tfrac{\deg \pi_m}{\deg \pi_n} 
(\beta_{m,n,\pi} \otimes \beta_{m,n,\pi^*})$.
That gives us direct systems of Hilbert spaces and the direct limits
\begin{equation}\label{lim0}
\cH_\pi = \varinjlim \{\cH_{\pi_n},\beta_{m,n,\pi}\} \text{ and } 
\cH_\pi\widehat{\otimes}\cH_\pi^* = \varinjlim \{(\cH_{\pi_n}\widehat{\otimes}
\cH_{\pi_n}^*),\zeta_{m,n,\pi}\},
\end{equation}
representation spaces for the irreducible unitary representations
$\pi = \varinjlim \pi_n$ of $G$ and $\pi\otimes\pi^*$ of $G\times G$.
\medskip

In order to sum the $\cH_{\pi_n}\widehat{\otimes}\cH_{\pi_n}^*$ and
the $\cH_\pi\widehat{\otimes}\cH_\pi^*$ as in (\ref{cent-decomp}), to form
$L^2(G_\chi) = \varinjlim L^2(G_{n,\chi})$, we need the direct integral
of (\ref{cent-decomp}) to be consistent with the rescaling isometries
of (\ref{incl-gen}).  Specifically we need
\begin{equation}\label{consistent0}
\begin{CD}
\int_{\widehat{G_{n,\chi}}} (\cH_{\pi_n}
        \widehat{\otimes}\cH_{\pi_n}^*) d_{n,\chi}(\pi_n)@>
\int_{\widehat{G_{n,\chi}}} \zeta_{m,n,\pi} d_{n,\chi}(\pi_n)>>
\int_{\widehat{G_{m,\chi}}} (\cH_{\pi_m}
        \widehat{\otimes}\cH_{\pi_m}^*) d_{m,\chi}(\pi_m)\\
@VV{\int_{\widehat{G_{n,\chi}}} \zeta_{m,n,\pi} d_{n,\chi}(\pi_n)}V  @VV{id}V\\
\int_{\widehat{G_{n,\chi}}} \zeta_{m,n,\pi}(\cH_{\pi_n}
        \widehat{\otimes}\cH_{\pi_n}^*) d_{n,\chi}(\pi_n)@>>>
\int_{\widehat{G_{m,\chi}}} (\cH_{\pi_m}
        \widehat{\otimes}\cH_{\pi_m}^*) d_{m,\chi}(\pi_m)
\end{CD}
\end{equation}
to commute for $m \geqq n$.
In the cases studied here, the $\widehat{G_{n,\chi}}$
are discrete (or even finite), so the $L^2(G_{n,\chi})$ are discrete direct
sums of irreducible 
representations.  Then there is no consistency problem with the rescaling
isometries: one simply takes their discrete direct sum.
\medskip

Once the $\cH_{\pi_n}\widehat{\otimes}\cH_{\pi_n}^*$ have been summed
as in (\ref{cent-decomp}), we need to control the integration over the
$\widehat{Z}$ in (\ref{planch--exp-finite-rel}) in order to pass to
the limit and form $L^2(G) = \int_{\widehat{Z}} L^2(G_\chi) d(\chi)$.
Here we need the conditions of (\ref{central-gen}), that every
$Z_n = Z$, and then we need that
\begin{equation}\label{abscont}
\text{for } m,n \ggg 0 \text{ the measures } d_m \text{ and } d_n
\text{ on } \widehat{Z} \text{ are mutually absolutely continuous.}
\end{equation}
In the cases studied here condition (\ref{abscont}) will be easy to verify.

\subsection{Construction of $\cA(G)$.}\label{sec2c}
Our hypothesis (\ref{req-sq-int}) of square integrability says that
Plancherel--almost--all $\cA(\pi_n) \subset \cH_{\pi_n}
\widehat{\otimes}\cH_{\pi_n}^*$, in fact form a dense subspace there.
We want the summation $L^2(G_{n,\chi}) = \int_{\widehat{G_{n,\chi}}} 
(\cH_{\pi_n} \widehat{\otimes}\cH_{\pi_n}^*) d_{n,\chi}(\pi_n)$
to restrict to a summation of the $\cA(\pi_n)$ and form a space
$\cA(G_{n,\chi})$ consisting of regular functions on $G_n$ that transform
by $\chi$.  In order to do this we must look at the detailed definition
of direct integral.
\medskip

\begin{definition} \label{def-l2-dirint}  {\em  Let $(Y, {\cM} ,\tau)$
be a measure space.  For each $y \in Y$ let $\cH_y$ be a separable 
Hilbert space.  Fix a family $\{s_\alpha\}_{\alpha \in A}$ of maps 
$Y \to \bigcup_{y \in Y} \cH_y$ such that
\begin{equation}\label{dirintmaps}
\begin{aligned}
\text{(i) } & s_\alpha(y) \in \cH_y \text{ a.e. }(Y, \cM ,\tau),
        \text{ for all } \alpha \in A, \\
\text{(ii) } & y \mapsto \langle s_\alpha(y), s_\beta(y) \rangle_{\cH_y}
        \text{ belongs to } L^1(Y, \tau), \text{ for all } \alpha ,\beta\in A,
        \text{ and }\\
\text{(iii) } & \cH_y \text{ is the closed span of }
        \{s_\alpha(y)\}_{\alpha \in A} \text{ a.e. } (Y,\tau).
\end{aligned} 
\end{equation}
Then the {\bf (Hilbert space) direct integral} defined by the measure
space $(Y, {\cM} ,\tau)$, the family $\{\cH_y \mid y \in Y\}$ of Hilbert
spaces, and the family $\{s_\alpha\}_{\alpha \in A}$ of maps, is the 
vector space
\begin{equation} \label{dirintspace}
\begin{aligned}
{\cH} = \int_Y \cH_y\,d\tau(y) : &\text{ all maps }
        s: Y \to \bigcup_{y \in Y} \cH_y \text{ such that } \\
\text{(i) } & s(y) \in \cH_y \text{ a.e. }(Y, \tau), \\
\text{(ii) } & y \mapsto \langle s(y), s_\alpha(y)\rangle_{\cH_y}
        \text{ is measurable, for each } \alpha \in A , \text{ and } \\
\text{(iii) } & y \mapsto \langle s(y), s_\alpha(y)\rangle_{\cH_y}
        \text{ belongs to } L^1(Y, \tau),
        \text{ for all } \alpha \in A \\
\text{with inner product} &\,\,\langle s , s' \rangle = \int_Y \langle s(y),
s'(y) \rangle_{\cH_y}\,d\tau(y).
\end{aligned} 
\end{equation}
}\hfill $\diamondsuit$ \end{definition}

The inner product of Definition \ref{def-l2-dirint} is well
defined, and the direct integral ${\cH}$ is a Hilbert space.
Our problem now is to find an appropriate family $\{s_\alpha\}_{\alpha \in A}$ 
of maps to $\bigcup_{\pi_n \in \widehat{G_n}} \cA(\pi_n)$ in order to put
the $\cA(\pi_n)$ together to make spaces $\cA(G_{n,\chi})$ and 
$\cA(G_n)$ of regular functions
on $G_n$ along the lines of a direct integral of Hilbert spaces.  
In other words we need the conditions that
\begin{equation}\label{fit-reg1}
\text{to form } L^2(G_{n,\chi}) = \int_{\widehat{G_{n,\chi}}} 
(\cH_{\pi_n} \widehat{\otimes}\cH_{\pi_n}^*) d_{n,\chi}(\pi_n)
	\text{ we may choose the } s_\alpha \text{ with } 
	s_\alpha(\pi_n) \in \cA(\pi_n)
\end{equation}
in order to form the 
$\cA(G_{n,\chi}) := \int_{\widehat{G_{n,\chi}}}\cA(\pi_n)d_{n,\chi}(\pi_n)$,
and then we need the conditions
{\small
\begin{equation}\label{fit-reg2}
\text{to form } L^2(G_n) = \int_{\widehat{Z}} L^2(G_{n,\chi}) d_n(\chi)
	\text{ we may assume } 
	s_\alpha\Bigl (\int_{\widehat{G_{n,\chi}}} (\pi_n\otimes\pi_n^*)
	d_{n,\chi}(\pi_n)\Bigr ) \in \cA(G_{n,\chi}).
\end{equation}
}
These conditions are automatic when the $G_n$ are compact. 
In Sections \ref{sec4} and \ref{sec5} we will verify them for the 
cases where $G_n$ is a connected simply connected 
nilpotent Lie group that has square integrable representations.
\medskip

Under the assumptions (\ref{consistent0}) and (\ref{abscont}) 
for forming our direct limits of Hilbert spaces, and (\ref{fit-reg1})
and (\ref{fit-reg2}) for restricting to regular functions before the
direct limits, we have the spaces
\begin{equation}\label{lim-reg}
\cA(G_n) := \int_{\widehat{Z}}\Bigl (\int_{\widehat{G_{n,\chi}}} \cA(\pi_n)
d_{n,\chi}(\pi_n)\Bigr )d(\chi) \text{ and }
\cA(G) := \varinjlim \cA(G_n).
\end{equation}
Here the $\cA(G_n)$ form a direct system using ordinary restriction of 
functions, and $\cA(G)$ is the direct limit of that system.  Each
$\cA(G_n)$ is a dense subspace of $L^2(G_n)$ but, because of rescaling
by the $\zeta_{m,n,\pi}$, we do not have $\cA(G)$ as a subspace of $L^2(G)$.
However there is a very useful relation which we now describe.

\subsection{Comparison of Direct Limits.}\label{sec2d}
In this section we will assume (\ref{consistent0}), (\ref{abscont}),
(\ref{fit-reg1}) and (\ref{fit-reg2}).  Then we have direct systems
$\{L^2(G_n)\} \to L^2(G)$ and $\{\cA(G_n)\} \to \cA(G)$.  Now define
$(G_n\times G_n)$--equivariant isometric inclusions
\begin{equation}\label{map-sys}
\eta_n: \cA(\pi_n) \to
\cH_{\pi_n}\widehat{\otimes}\cH_{\pi_n}^* \text{ by }
f_{u,v,n} \mapsto \sqrt{\deg \pi_n}\, f_{u,v,n}.
\end{equation}
\begin{proposition}\label{comparison}
The maps $\eta_n$ of {\rm (\ref{map-sys})} satisfy
$\eta_m\circ\zeta_{m,n,\pi}(f_{u,v,n}) = 
\zeta_{m,n,\pi}\circ\eta_n(f_{u,v,n})$ and
send the direct system $\{\cA(G_n)\}$
into the direct system $\{L^2(G_n)\}$.  That map of direct systems defines a 
$(G\times G)$--equivariant injection 
$$
\eta: \cA(G) \to L^2(G)
$$ 
with dense image.  In particular $\eta$ defines a pre Hilbert space
structure on $\cA(G)$ with completion isometric to $L^2(G)$.
\end{proposition}
\noindent {\bf Proof.}  We extract the result from \cite[Appendix A]{NRW}.
Let $\cC$ denote the category of topological vector spaces and continuous
linear maps.  We compute $\eta_m\circ\zeta_{m,n,\pi}(f_{u,v,n}) = \eta_n(f_{u,v,n})$,
as asserted, and use (\ref{consistent0}) and (\ref{abscont}) to
patch the $\eta_n$ together to form maps 
$\widetilde{\eta_n}: \cA(G_n) \to L^2(G_n)$.
The space $L^2(G_n)$ carries its usual topology and we give $\cA(G_n)$
the subspace topology.  Now view $\{\cA(G_n)\}$ and $\{L^2(G_n)\}$
as direct systems in the category $\cC$.  The $\widetilde{\eta_n}$ define a
morphism $\{\cA(G_n)\} \to \{L^2(G_n)\}$ of direct systems in $\cC$,
so the universal property of direct limits gives a morphism
$\widetilde{\eta}:\varinjlim_\cC \{\cA(G_n)\} \to \varinjlim_\cC \{L^2(G_n)\}$
of the direct limits in $\cC$.  Note that $\widetilde{\eta_n}$ is injective
and $(G_n \times G_n)$--equivariant with dense image.  It follows that 
$\widetilde{\eta}$ is injective and $(G\times G)$--equivariant with dense image.
\medskip

The topological vector space direct limit $\varinjlim_\cC \{L^2(G_n)\}$
has a pre Hilbert space structure given by the Hilbert space structures
on the $L^2(G_n)$, and the Hilbert space direct limit 
$L^2(G) = \varinjlim \{L^2(G_n)\}$ is its Hilbert space completion. Thus
$\widetilde{\eta}:\varinjlim_\cC \{\cA(G_n)\} \to \varinjlim_\cC \{L^2(G_n)\}$
is in fact a continuous linear map, injective and $(G\times G)$--equivariant,
from $\varinjlim_\cC \{\cA(G_n)\}$ to $L^2(G_n)$.  Further, our original
$\cA(G)$ is obtained from $\varinjlim_\cC \{\cA(G_n)\}$ by forgetting the
topology.  Thus in particular $\widetilde{\eta}$ maps $\cA(G)$ into 
$L^2(G)$, and it is a $(G\times G)$--equivariant injection onto a dense 
subspace.
\hfill $\square$

\subsection{Homogeneous Spaces.}\label{sec2e}
We now consider direct limit homogeneous spaces $G/K = \varinjlim G_n/K_n$.
Specifically, we require that $G = \varinjlim G_n$, that the $K_n$ are compact
subgroups of the $G_n$ such that $K_n = G_n\cap K_m$ for $m \geqq n$, and
that $K = \varinjlim K_n$.  We also assume (\ref{consistent0}), (\ref{abscont}),
(\ref{fit-reg1}) and (\ref{fit-reg2}) so that Proposition \ref{comparison}
applies to $\{G_n\}$ and $G$.
\medskip

Since $K_n$ is compact, $L^2(G_n/K_n) = L^2(G_n)^{K_n}$, the subspace
consisting of the right--$K_n$--invariant $L^2$ functions on $G_n$.  Of
course we can define $L^2(G/K) = L^2(G)^K$, but for comparison with
$\cA(G/K)$ we want $L^2(G/K)$ to be of the form $\varinjlim L^2(G_n/K_n)$.
Since $K_n$ is compact, (\ref{planch--exp-finite}) gives us 
\begin{equation}\label{rightinv}
L^2(G_n/K_n) = L^2(G_n)^{K_n} = \int_{\widehat{G_n}} 
\left (\cH_{\pi_n} \widehat{\otimes} (\cH_{\pi_n}^*)^{K_n}\right ) d(\pi_n).
\end{equation}
we view the direct limit maps $\beta_{m,n,\pi}: \cH_{\pi_n} \to \cH_{\pi_m}$
as isometric inclusions and
arrange the $L^2(G_n/K_n)$ into a direct system by requiring that
\begin{equation}\label{k-fixed-consistent}
\text{orthogonal projection } 
p_{m,n,\pi}: \cH_{\pi_m} \to \cH_{\pi_n}
\text{ defines a bijective map } \cH_{\pi_m}^{K_m} \to \cH_{\pi_n}^{K_n}.
\end{equation}
$\text{Then every }u_m \in \cH_{\pi_m}^{K_m} \text{ has unique expression }
u_m = u_n + x \text{ with } u_n \in \cH_{\pi_n}^{K_n}
\text{ and } x \perp \beta_{m,n,\pi}(\cH_{\pi_n})$.
As $\pi_n$ is irreducible and $||u_m||^2 = ||u_n||^2 + ||x||^2$
we have a constant $c = c_{m,n,\pi}$ such that
$||u_n|| = c||u_m||$.  Here $0 < c \leqq 1$ and $u_n = p_{m,n,\pi}(u_m)$.
Thus we have an isometry
\begin{equation}\label{scale-invars}
\alpha_{m,n,\pi}: \cH_{\pi_n}^{K_n} \cong \cH_{\pi_m}^{K_m} \text{ by }
\alpha_{m,n,\pi}(u_n) = c_{m,n,\pi}u_m
\text{ where } u_n = p_{m,n,\pi}(u_m).
\end{equation}
Now we have $G_n$--equivariant isometric injections
\begin{equation}\label{scale-invars1}
\widetilde{\zeta}_{m,n,\pi}: \cH_{\pi_n}\widehat{\otimes}(\cH_{\pi_n}^*)^{K_n} 
\to \cH_{\pi_m}\widehat{\otimes}(\cH_{\pi_m}^*)^{K_m} \text{ defined by }
f_{u,v_n,n} \mapsto 
c_{m,n,\pi}\left ( \tfrac{\deg \pi_m}{\deg \pi_n}\right )^{1/2}f_{u,v_m,m}.
\end{equation}
Then $\left \{\left (\cH_{\pi_n} \widehat{\otimes}(\cH_{\pi_n}^*)^{K_n}\right ),
\widetilde{\zeta}_{m,n,\pi}\right\}$ is isomorphic to a subsystem 
of the Hilbert space direct system
$\left \{\left (\cH_{\pi_n} \widehat{\otimes}\cH_{\pi_n}^*\right ),
\zeta_{m,n,\pi}\right \}$ of (\ref{lim0}).  As a result we have
\begin{equation}\label{rightinv1}
\cH_\pi \widehat{\otimes} (\cH_\pi^*)^K = \varinjlim
\left \{\left (\cH_{\pi_n} \widehat{\otimes}(\cH_{\pi_n}^*)^{K_n}\right ),
\widetilde{\zeta}_{m,n,\pi}\right \} 
\end{equation}
in the Hilbert space category,
and they fit together under the direct integral (\ref{rightinv})
to give us
\begin{equation}\label{rightinv2}
L^2(G/K) := L^2(G)^K = 
\varinjlim \left \{L^2(G_n/K_n), \widetilde{\zeta}_{m,n}\right \}.
\end{equation}

From (\ref{lim-coef}) and (\ref{k-fixed-consistent}) we have 
\begin{equation}\label{quo-lim-coef}
\begin{aligned}
&\cA(\pi_n)^{K_n} = \{\text{finite linear combinations of the } f_{u,v_n,n}
  \text{ where } u \in \cH_{\pi_n} \text{ and } v_n \in \cH_{\pi_n}^{K_n}\} 
	\text{ and }\\
&\nu_{m,n,\pi}:\cA(\pi_n)^{K_n} \hookrightarrow \cA(\pi_m)^{K_m} 
	\text{ by }
	f_{u,v_n,n} \mapsto f_{u,v_m,m} \text{ where projection } 
	p_{m,n,\pi}(v_m) = v_n.
\end{aligned}
\end{equation}
\begin{lemma}\label{res-lemma}
If $f \in \cA(\pi_n)^{K_n}$ then $\nu_{m,n,\pi}(f)|_{G_n/K_n} = f$.
\end{lemma}
\noindent {\bf Proof.} Let $u \in \cH_{\pi_n}$ and $v_n \in \cH_{\pi_n}^{K_n}$. 
Then $v_m = v_n + x \in \cH_{\pi_m}^{K_m}$ with $x \perp \cH_{\pi_n}$.
Given $g \in G_n$ compute
$\left (\nu_{m,n,\pi}(f_{u,v_n,n})\right )(g) 
= f_{u,v_m,m}(g) 
= \langle u, \pi_m(g)(v_n + x)\rangle 
= \langle \pi_m(g^{-1})u, v_n + x)\rangle
= \langle \pi_m(g^{-1})u, v_n \rangle
= \langle u, \pi_m(g)(v_n)\rangle 
= \langle u, \pi_n(g)(v_n)\rangle = f_{u,v_n,n}(g)$.
\hfill $\square$

Equation \ref{quo-lim-coef} defines a direct system 
$\{\cA(\pi_n)^{K_n},\nu_{m,n,\pi}\}$ and Lemma \ref{res-lemma} shows that
its maps are inverse to restriction.  Thus the direct limit
\begin{equation}\label{quo-lim-coef1}
        \cA(\pi)^K := \varinjlim \{\cA(\pi_n)^{K_n}, \nu_{m,n,\pi}\}
\end{equation}
sits naturally as a $G$--submodule of the inverse limit
$\varprojlim \{\cA(\pi_n)^{K_n},\text{ restriction}\}$.
Now, from Proposition \ref{comparison}, we can take right $K$--invariants
as follows.

Our assumptions (\ref{consistent0}) and (\ref{abscont})
for forming our direct limits of Hilbert spaces, and (\ref{fit-reg1})
and (\ref{fit-reg2}) for restricting to regular functions before the
direct limits, carry (\ref{lim-reg}) over to $K_n$-- and $K$--invariant
regular functions as follows. 
\begin{equation}\label{quo-lim-reg}
\begin{aligned}
\cA(G_n/K_n) := &\cA(G_n)^{K_n}
= \int_{\widehat{Z}}\Bigl (\int_{\widehat{G_{n,\chi}}} 
	\cA(\pi_n)^{K_n} d_{n,\chi}(\pi_n)\Bigr )d(\chi)\\
&\text{ and }
\cA(G/K) := \cA(G)^K = 
   \varinjlim \{\cA(G_n/K_n), \nu_{m,n}\}.
\end{aligned}
\end{equation}
In view of Lemma \ref{res-lemma} the $\cA(G_n/K_n)$ form a direct system 
whose maps are inverse to restriction of functions 
and $\cA(G/K)$ is the direct limit of that system.  Each
$\cA(G_n/K_n)$ is a dense subspace of $L^2(G_n/K_n)$ but, because the
$\nu_{m,n}$ are not isometric, we do not have $\cA(G/K)$ sitting naturally 
as a subspace of $L^2(G/K)$.
As in the group level setting we can manage this in a somewhat abstract way.
The $G_n$--equivariant isometric inclusions restrict to
\begin{equation}\label{rel-map-sys}
\widetilde{\eta}_{n,\pi}: \cA(\pi_n)^{K_n} \to
\cH_{\pi_n}\widehat{\otimes}(\cH_{\pi_n}^*)^{K_n} \text{ by }
f_{u,v,n} \mapsto c_{n,1,\pi}\sqrt{\deg \pi_n}\, f_{u,v,n}.
\end{equation}

\begin{proposition}\label{quo-comparison}
The maps $\widetilde{\zeta}_{m,n,\pi}$ of {\rm (\ref{scale-invars1})},
$\nu_{m,n,\pi}$ of {\rm (\ref{quo-lim-coef})} and
$\widetilde{\eta}_{n,\pi}$ of {\rm (\ref{rel-map-sys})}
satisfy $$(\widetilde{\eta}_{m,\pi}\circ \nu_{m,n,\pi})(f_{u,v,n}) = 
(\widetilde{\zeta}_{m,n,\pi}\circ\widetilde{\eta}_{n,\pi})(f_{u,v,n})$$
for $f_{u,v,n} \in \cA(\pi_n)^{K_n}$.  Thus they inject
the direct system $\{\cA(G_n)^{K_n}, \nu_{m,n}\}$
into the direct system $\{L^2(G_n)^{K_n},\widetilde{\zeta}_{m,n}\}$.  
That map of direct systems defines a $G$--equivariant injection 
$$
\widetilde{\eta}: \cA(G/K) \to L^2(G/K)
$$
with dense image.  In particular $\widetilde{\eta}$ defines a pre Hilbert space
structure on $\cA(G/K)$ with completion isometric to $L^2(G/K)$.
\end{proposition}
\noindent {\bf Proof.}  The point is that $c_{m,n,\pi} 
= c_{m,m-1,\pi}\cdot c_{m-1,m-2,\pi}\cdot \ldots\cdot c_{n+1,n,\pi}$; so
$c_{m,n,\pi} = c_{m,1,\pi}/c_{n,1,\pi}$.  Now compute
$$
\begin{aligned}
(\widetilde{\zeta}_{m,n,\pi}\circ\widetilde{\eta}_{n,\pi})(f_{u,v_n,n})
   &= \widetilde{\zeta}_{m,n,\pi}(c_{n,1,\pi}\sqrt{\deg\pi_n}\, f_{u,v_n,n})\\
&= c_{m,n,\pi}\left ( \tfrac{\deg\pi_m}{\deg\pi_n}\right )^{1/2}
		c_{n,1,\pi}\sqrt{\deg\pi_n}\, f_{u,v_m,m} \\
&= c_{m,1,\pi}\sqrt{\deg\pi_m}\, f_{u,v_m,m} 
= (\widetilde{\eta}_{m,\pi}\circ\nu_{m,n,\pi})(f_{u,v_n,n}).
\end{aligned}
$$
Our assertions follow as in the proof of Proposition \ref{comparison}.
\hfill $\square$

\subsection{Commutative Spaces.}\label{sec2f}
Suppose that the $G_n/K_n$ are commutative spaces, i.e. that the
$(G_n,K_n)$ are Gelfand pairs.  If $\pi_n \in \widehat{G_n}$
then $\dim \cH_{\pi_n}^{K_n} \leqq 1$, in other words either
$\cH_{\pi_n}^{K_n} = 0$ or there is a unique (up to scale) unit vector
$v_n \in \cH_{\pi_n}^{K_n}$.  In the latter case 
(\ref{k-fixed-consistent}) simply asserts that $v_m$ cannot be orthogonal
to $\beta_{m,n,\pi}(v_n)$, so that Proposition \ref{quo-comparison}
is valid.  Here there is a problem: given $\cH_{\pi_n}^{K_n} \ne 0$ and 
$m \geqq n$ we must have $\cH_{\pi_m}^{K_m} \ne 0$.  
In the remainder of this paper we study two situations in which this
holds: when the $G_n/K_n$ are compact symmetric spaces, and in
many cases where the $G_n/K_n$ are commutative nilmanifolds.

\section{Limits of Compact Lie Groups and Compact Symmetric Spaces}
\label{sec3}
\setcounter{equation}{0}

Suppose that $G_n$ is a compact topological group.  Then the Peter Weyl 
Theorem says that $L^2(G_n)$ is the Hilbert space direct sum of the
spaces $\cA(\pi_n)$ of coefficients, $\pi_n \in \widehat{G_n}$.  In
particular one has a dense subspace of $L^2(G_n)$ given by 
the algebraic direct sum
\begin{equation}\label{PW}
\cA(G_n) = {\sum}_{\pi_n \in \widehat{G_n}}\, \cA(\pi_n).
\end{equation}
When $G_n$ is a compact Lie group, the spectrum of the ring $\cA(G_n)$
is a linear algebraic group, the {\it associated algebraic group} of $G_n$
for which $\cA(G_n)$ is the ring of regular functions.  This essentially
is Tannaka's Theorem; see Chevalley \cite[Chapter 6]{Ch}.
Thus $\cA(G_n)$ is the ring of regular functions on $G_n$
as well as on its associated algebraic group.  Note that the 
associated algebraic group is the complexification $(G_n)_\C$.

\subsection{Parabolic Direct Limits of Groups.}\label{sec3a}
We make the standing assumption for this section and the next that
\begin{equation}\label{req-parabolic}
\begin{aligned}
&\text{the } G_n \text{ are compact connected Lie groups and } \\
&\{G_n, \varphi_{m,n}\} \text{ is a strict parabolic direct system.}
\end{aligned}
\end{equation}
Here {\it parabolic} means that the semisimple part of $\varphi_{m,n}(G_n)$
is the semisimple part of the centralizer of a toral subgroup of $G_m$,
in other words that $\varphi_{m,n}([G_n,G_n])_\C$ is the semisimple
component of a parabolic subgroup of $(G_m)_\C$.  Thus we have Cartan
subalgebras $\gh_n \subset \gg_n$ such that $d\varphi_{m,n}(\gh_n)
\subset \gh_m$ and simple root systems $\Psi_n \subset i\gh_n^*$ such
that
\begin{equation}\label{res-simple-roots0}
\text{if } \psi \in \Psi_n \text{ then there is a unique } \psi' \in \Psi_m
\text{ such that } d\varphi_{m,n}^*(\psi') = \psi.
\end{equation}
Now we can enumerate the simple root systems as
\begin{equation}\label{res-simple-roots1}
\Psi_n = \{\psi_{n,1}, \dots , \psi_{n,\ell_n}\} \, ,
\ell_n = \rank \gg_n \text{, and } \psi_{n,j} = d\varphi_{m,n}^*(\psi_{m,j})
\text{ for } 1 \leqq j \leqq \ell_n \text{ and } m \geqq n.
\end{equation}

Dually we have enumerations of the systems of fundamental weights as
\begin{equation}\label{res-simple-wts1}
\Xi_n = \{\xi_{n,1}, \dots , \xi_{n,\ell_n}\} \, ,
\ell_n = \rank \gg_n \text{, and } \xi_{n,j} = d\varphi_{m,n}^*(\xi_{m,j})
\text{ for } 1 \leqq j \leqq \ell_n \text{ and } m \geqq n.
\end{equation}
Now the correspondence 
$$
k_1\xi_{n,1} + \dots + k_{\ell_n}\xi_{n,\ell_n} \mapsto
k_1\xi_{m,1} + \dots + k_{\ell_m}\xi_{n,\ell_n}
$$
sends dominant integral weights to dominant integral weights.  Fix an
index $n_0$ and a dominant integral weight $\lambda_{n_0}$ for
$(\gg_{n_0}, \gh_{n_0}, \Psi_{n_0})$.  For $n \geqq n_0$ we now have
a dominant integral weight $\lambda_n$ for $(\gg_n,\gh_n,\Psi_n)$
and a corresponding direct system of irreducible representations of the
$G_n$ given by
\begin{equation}\label{dirsys-reps}
\begin{aligned}
\{\pi_{n,\lambda}\} \text{ where } 
\lambda_{n_0} &= k_1\xi_{n_0,1} + \dots + k_{\ell_{n_0}}\xi_{n_0,\ell_{n_0}}
\text{ and } \\
&\pi_{n,\lambda} \in \widehat{G_n} \text{ has highest weight }
\lambda_n = k_1\xi_{n,1} + \dots + k_{\ell_{n_0}}\xi_{n,\ell_{n_0}}.
\end{aligned}
\end{equation}

Choose unit highest weight vectors $v_{n,\lambda} \in \cH_{\pi_{n,\lambda}}$.
We embed representation spaces 
\begin{equation}\label{zetaprime}
\beta_{m,n,\lambda} = \beta_{m,n,\pi_\lambda}: 
\cH_{\pi_{n,\lambda}} \hookrightarrow
\cH_{\pi_{m,\lambda}} \text{ by } X(v_{n,\lambda}) \mapsto X(v_{m,\lambda})
\text{ for } X \in \cU(\gg_n).
\end{equation}
Then we have the direct limit representation and its representation space:
\begin{equation}\label{limrep}
\pi_\lambda = \varinjlim \pi_{n,\lambda} \text{ unitary representation of }
G = \varinjlim G_n \text{ on } 
\cH_{\pi_\lambda} = \varinjlim \cH_{\pi_{n,\lambda}}\, .
\end{equation}
As in (\ref{incl-gen}) we now have isometric injections
\begin{equation}\label{normalize-inj-cpt}
\zeta_{m,n}: \cH_{\pi_{n,\lambda}}\widehat{\otimes}\cH_{\pi_{n,\lambda}}^* \to
\cH_{\pi_{m,\lambda}}\widehat{\otimes}\cH_{\pi_{m,\lambda}}^* \text{ defined by }
f_{u,v,n} \mapsto \left ( \tfrac{\deg \pi_{m,\lambda}}{\deg \pi_{n,\lambda}}
\right )^{1/2}f_{u,v,m}.
\end{equation}
That gives us direct systems of Hilbert spaces and the direct limits
\begin{equation}\label{iso-dir-lim-cpt}
\cH_{\pi_\lambda}\widehat{\otimes}\cH_{\pi_\lambda}^* = 
\varinjlim \{(\cH_{\pi_{n,\lambda}}\widehat{\otimes} 
\cH_{\pi_{n,\lambda}}^*),\zeta_{m,n}\}. 
\end{equation}
They are the representation spaces for the irreducible unitary representations
$\pi_\lambda$ of $G$ and $\pi_\lambda\otimes\pi_\lambda^*$ of $G\times G$.
The Peter--Weyl Theorem for Parabolic Direct Limits \cite[Theorem 4.3]{W5}
exhibits the Hilbert space $L^2(G) := \varinjlim \{L^2(G_n),\zeta_{m,n}\}$ 
as the orthogonal
direct sum of the $\cH_{\pi_\lambda}\widehat{\otimes}\cH_{\pi_\lambda}^*$,
and shows that the left/right regular representation of $G\times G$ 
on $L^2(G)$ the multiplicity--free 
discrete direct sum of irreducible representations
$\pi_\lambda \otimes \pi_\lambda^*$ of highest weights $(\lambda,\lambda^*)$.
\medskip

The direct integral conditions (\ref{consistent0}) and (\ref{abscont})
are automatic here because the integrals are direct sums.  Thus, as
$\lambda$ varies, the (\ref{normalize-inj-cpt}) and their limits 
(\ref{iso-dir-lim-cpt}) sum to give us
\begin{equation}\label{n-normalize-inj-cpt}
\zeta_{m,n}: L^2(G_n) \to L^2(G_m) \text{ defined by }
f_{u,v,n} \mapsto \left ( \tfrac{\deg \pi_{m,\lambda}}{\deg \pi_{n,\lambda}}
\right )^{1/2}f_{u,v,m} \text{ for } u, v \in \cH_{\pi_{n,\lambda}}
\end{equation}
and the limit Hilbert space
\begin{equation}\label{i-iso-dir-lim-cpt}
L^2(G) := 
\varinjlim \{ L^2(G_n), \zeta_{m,n}\}.
\end{equation}

\subsection{The Ring of Regular Functions for Parabolic Direct Limits.}
\label{sec3b}
As noted earlier, the direct integral conditions
(\ref{consistent0}), (\ref{abscont}), (\ref{fit-reg1}) and
(\ref{fit-reg2}) are automatic in the compact group setting.  Thus 
Proposition \ref{comparison} applies, showing that the
$\cA(G_n)$ and $\cA(G)$ follow the same pattern as in the Peter--Weyl 
Theorem for Parabolic Direct Limits, and the resulting map
$\cA(G) \to L^2(G)$ defines both a pre Hilbert space structure on
$\cA(G)$ and an interpretation of the elements of $L^2(G)$ as the Hilbert
space completion of $\cA(G)$.

\subsection{Parabolic Direct Limits of Compact Symmetric Spaces}\label{sec3c}

Fix a parabolic direct system of compact connected Lie groups $G_n$ and 
subgroups
$K_n$ such that each $(G_n,K_n)$ is an irreducible riemannian symmetric pair.
Suppose that the corresponding compact symmetric spaces $M_n = G_n/K_n$ are
simply connected.  Up to re--numbering and passage to a common
cofinal subsequence the only possibilities are
{\small
\begin{equation}\label{symmetric-case-class}
\begin{tabular}{|c|l|l|c|c|} \hline
\multicolumn{5}{| c |}
{Compact Irreducible Riemannian Symmetric Spaces $M_n = G_n/K_n$} \\
\hline \hline
\multicolumn{1}{|c}{} &
        \multicolumn{1}{c}{$G_n$} &
        \multicolumn{1}{|c}{$K_n$} &
        \multicolumn{1}{|c}{Rank$M_n$} &
        \multicolumn{1}{|c|}{Dim$M_n$} \\ \hline \hline
$1$ & $SU(n)\times SU(n)$ & diagonal $SU(n)$ & $n-1$ & $n^2-1$ \\ \hline
$2$ & $Spin(2n+1)\times Spin(2n+1)$ & diagonal $Spin(2n+1)$ &
        $n$ & $2n^2+n$ \\ \hline
$3$ & $Spin(2n)\times Spin(2n)$ & diagonal $Spin(2n)$ &
        $n$ & $2n^2-n$ \\ \hline
$4$ & $Sp(n)\times Sp(n)$ & diagonal $Sp(n)$ & $n$ & $2n^2+n$ \\ \hline
$5$ & $SU(p+q), \ p = p_n, q = q_n$ & $S(U(p)\times U(q))$ &
        $\min(p,q)$ & $2pq$ \\ \hline
$6$ & $SU(n)$ & $SO(n)$ & $n-1$ & $\frac{(n-1)(n+2)}{2}$ \\ \hline
$7$ & $SU(2n)$ & $Sp(n)$ & $n-1$ & $2n^2-n-1$  \\ \hline
$8$ & $SO(p+q), \ p = p_n, q = q_n$ & $SO(p) \times SO(q)$ &
        $\min(p,q)$ & $pq$  \\ \hline
$9$ & $SO(2n)$ & $U(n)$ & $[\frac{n}{2}]$ & $n(n-1)$ \\ \hline
$10$ & $Sp(p+q), \ p = p_n, q = q_n$ & $Sp(p) \times Sp(q)$ &
        $\min(p,q)$ & $4pq$  \\ \hline
$11$ & $Sp(n)$ & $U(n)$ & $n$ & $n(n+1)$  \\ \hline
\end{tabular}
\end{equation}
}
These are the cases where the $G_n$ form a parabolic direct system.
Now we set about carrying the results of Sections \ref{sec3a} and 
\ref{sec3b} from $G = \varinjlim G_n$ to $M = \varinjlim M_n$ for 
the systems of Table \ref{symmetric-case-class}.

\subsection{Square Integrable Functions and Regular Functions.}\label{sec3d}

Recall the decomposition $\gg_n = \gk_n + \gs_n$ under the symmetry
$\theta$ of $M_n$.  Here $K_n$ is the identity component of
the fixed point set $G_n^\theta$, $\gk_n$ is the $+1$ eigenspace of $d\theta$,
and $\gs_n$ is the $-1$ eigenspace.  We assume the alignments
$\theta_n = \theta_m|_{G_n}$ so $\gk_n \subset \gk_m$ and $\gs_n \subset \gs_m$.
We recursively choose maximal abelian subspaces $\ga_n \subset \gs_n$ with
$\ga_n \subset \ga_m$ and define $\gm_n$ to be the centralizer of $\ga_n$
in $\gk_n$.  For the systems of Table \ref{symmetric-case-class} we have
$\gm_n \subset \gm_m$, and we recursively choose Cartan subalgebras
$\gt_n \subset \gm_n$ with that $\gt_n \subset \gt_m$.  Then the 
$\gh_n := \gt_n + \ga_n$ are $\theta$--stable Cartan subalgebras of the $\gg_n$.
\medskip

The restricted root systems
$$
\Delta(\gg_n,\ga_n) = \{\alpha|_{\ga_n} \mid \alpha \in \Delta(\gg_n,\gh_n)
\text{ and } \alpha|_{\ga_n} \ne 0\}
$$
have consistent root orders
$$
\text{if } \alpha \in \Delta^+(\gg_m,\ga_m) \text{ and } \alpha|_{\ga_n} \ne 0
\text{ then } \alpha|_{\ga_n} \in \Delta^+(\gg_n,\ga_n),
$$
and similarly we have consistent root orders on the $\Delta(\gm_n,\gt_n)$.
Together they define consistent positive root systems
$$
\Delta^+(\gg_n,\gh_n) = \{\alpha \in \Delta(\gg_n,\gh_n) \mid
\text{ either } \alpha|_{\ga_n} \in \Delta^+(\gg_n,\ga_n), \text{ or }
\alpha|_{\ga_n} = 0 \text{ and } \alpha|_{\gt_n} \in \Delta^+(\gm_n,\gt_n)\}.
$$
Using the parabolic condition we have simple restricted root systems
$$
\Psi_n = \Psi(\gg_n,\ga_n) \text{ such that if } \psi \in \Psi_n
\text{ then there is a unique } \psi' \in \Psi_m \text{ such that }
\psi = \psi'|_{\ga_n}.
$$
Now as in (\ref{res-simple-roots1}) we enumerate
\begin{equation}\label{res-simple-res-roots1}
\Psi_n = \{\psi_{n,1}, \dots , \psi_{n,\ell_n}\} \, ,
\ell_n = \dim \ga_n \text{, and } \psi_{n,j} = \psi_{m,j}|_{\ga_n}
\text{ for } 1 \leqq j \leqq \ell_n \text{ and } m \geqq n.
\end{equation}

Using the root orders just described, the 
Cartan--Helgason Theorem says that the irreducible representation
$\pi_\lambda$ of $\gg_n$ of highest weight $\lambda$ gives a summand of
the representation of $G_n$ on $L^2(M_n)$ if and only if (i)
$\lambda|_{\gt_n} = 0$, so we may view $\lambda$ as an element of
$i\ga_n^*$, and (ii) if $\alpha \in \Delta^+(\gg_n , \ga_n)$ then
$\tfrac{\langle \lambda, \alpha \rangle} {\langle \alpha, \alpha \rangle}$
is an integer $\geqq 0$.  
\medskip

Condition (i) of the Cartan--Helgason Theorem persists under
restriction $\lambda \mapsto \lambda|_{\gh_{n-1}}$ because
$\gt_{n-1} \subset \gt_n$.  Given (i), condition (ii) says that
$\tfrac{1}{2}\lambda$ belongs to the weight lattice of $\gg_n$, so its
restriction to $\gh_{n-1}$ exponentiates to a well defined function on
the corresponding maximal torus of $G_{n-1}$ and thus belongs to the
weight lattice of $\gg_{n-1}$.  Given condition (i), our arrangements
$\ga_n \subset \ga_m$ and $\gt_n \subset \gt_m$ ensure that condition (ii)
persists under restrictions $\lambda \mapsto \lambda|_{\gh_{n-1}}$.
With this in mind, we
define linear functionals $\xi_{n,j}\in i\ga_n^*$ by
\begin{equation}\label{res-simple-res-weights1}
\tfrac{\langle \xi_{n,i},\psi_{n,j} \rangle}
{\langle \psi_{n,j},\psi_{n,j}\rangle}=\delta_{i,j} \text{ for }
1 \leqq j \leqq r_n \text{ --- except that }
\tfrac{\langle \xi_{n,i},\psi_{n,i} \rangle}
{\langle \psi_{n,i},\psi_{n,i}\rangle} = 2 \text{ if }
2\psi_{n,i} \in \Delta(\gg_n,\ga_n).
\end{equation}
The weights $\xi_{n,j}$ are the
\textit{class 1 fundamental highest weights} for $(\gg_n,\gk_n)$
Now as in Section \ref{sec3a} the correspondence
$$
k_1\xi_{n,1} + \dots + k_{\ell_n}\xi_{n,\ell_n} \mapsto
k_1\xi_{m,1} + \dots + k_{\ell_m}\xi_{n,\ell_n}
$$
sends class $1$ dominant integral weights to class $1$ dominant integral 
weights and defines direct systems of irreducible class $1$ representations by
\begin{equation} \label{dirsys-class1-reps}
\begin{aligned}
\{\pi_{n,\lambda}\} \text{ where }
\lambda_{n_0} &= k_1\xi_{n_0,1} + \dots + k_{\ell_{n_0}}\xi_{n_0,\ell_{n_0}}
\text{ and } \\
&\pi_{n,\lambda} \in \widehat{G_n} \text{ has highest weight }
\lambda_n = k_1\xi_{n,1} + \dots + k_{\ell_{n_0}}\xi_{n,\ell_{n_0}}.
\end{aligned}
\end{equation}
As before, the corresponding embeddings of representation spaces 
$\beta_{m,n,\lambda}: \cH_{\pi_{n,\lambda}} \to \cH_{\pi_{n,\lambda}}$
are given by choices of highest weight unit vectors $v_{n,\lambda}
\in \cH_{\pi_{n,\lambda}}$ and by
$X(v_{n,\lambda}) \mapsto  X(v_{m,\lambda})$ for $X \in \cU(\gg_n)$.
Again we have the direct limit representations and its representation space
\begin{equation}\label{dirlim-class1}
\pi_\lambda = \varinjlim \pi_{n,\lambda} \text{ unitary representation of }
G = \varinjlim G_n \text{ on } \cH_{\pi_\lambda} = \varinjlim
\{\cH_{\pi_{n,\lambda}}, \beta_{m,n,\lambda}\}.
\end{equation}
Note that $\pi_\lambda$ is irreducible, has highest weight $\lambda$, and 
has a highest weight unit vector $v_\lambda = \varinjlim v_{n,\lambda}$ 
that is invariant under the action of $K = \varinjlim K_n$.  Given
$\{\lambda_n\}$ as in (\ref{dirsys-class1-reps}) we write
$\cH_{n,\lambda}$ for $\cH_{\pi_{n,\lambda}}$.

\begin{lemma}\label{invariants}
For $\{\lambda_n\}$ as in {\rm (\ref{dirsys-class1-reps})} and each 
$n \geqq n_0$
let $w_{n,\lambda}$ be a unit vector in $\cH_{n,\lambda}^{K_n}$.  Then
orthogonal projection $p_{m,n,\lambda}: \cH_{m,\lambda} \to \cH_{n,\lambda}$
maps $w_{m,\lambda}$ to a nonzero multiple $c_{m,n,\lambda}w_{n,\lambda}$
of $w_{n,\lambda}$.
\end{lemma}
\noindent {\bf Proof.} 
By $G_n$--equivariance of the projection, 
$p_{m,n,\lambda}(w_{m,\lambda})$ is $K_n$--invariant, hence is a multiple of  
$w_{n,\lambda}$.  Suppose that the multiple is zero, i.e. that
$w_{n,\lambda} \perp w_{m,\lambda}$.  Then $w_{m,\lambda} \perp
\cH_{n,\lambda}$, in particular $w_{m,\lambda}$ is orthogonal to the
highest weight $\lambda$ vector $v_{n,\lambda}$.
As $v_{m,\lambda} = v_{n,\lambda}$, now $w_{m,\lambda} \perp v_{m,\lambda}$.
But the argument of \cite[Proposition 2.2]{Li} shows
that $w_{m,\lambda}$ cannot be orthogonal to $v_{m,\lambda}$.  We conclude
$p_{m,n,\lambda}(w_{m,\lambda}) = c w_{n,\lambda}$ for some nonzero
constant $c$.
\hfill $\square$
\medskip

Denote $\lambda^* = -w_0(\lambda)$ where $w_0$ is the longest element of
the Weyl group.  Then $\lambda^*$ is the highest weight of the dual 
$\pi_{n,\lambda}^*$, and
$\cH_{\pi_{n,\lambda}}^* = \cH_{\pi_{n,\lambda^*}}$.
Given $\pi_{n,\lambda} \in \widehat{G_{n,K}}$ we have
$\pi_{n,\lambda^*} \in \widehat{G_{n,K}}$ and we may suppose that its
$K_n$--fixed unit vector $w_{\lambda^*}$ is normalized so that
the pairing $(w_\lambda, w_{\lambda^*}) = 1$.  Then the
space of right--$K_n$--invariant matrix coefficient functions corresponding
to $\pi_{n,\lambda}$ is (as a module for the left translation action of $G_n$)
\begin{equation}\label{coefs-on-symm}
\left ( \cH_{\pi_{n,\lambda}} \otimes \cH_{\pi_{n,\lambda}}^* \right )^{K_n}
\ = \ \cH_{\pi_{n,\lambda}} \otimes w_{\lambda^*}\C
\ = \ \{f_{u,w_\lambda,n} \mid u \in \cH_{\pi_{n,\lambda}}\}
\ \cong \ \cH_{\pi_{n,\lambda}}.
\end{equation}

Lemma \ref{invariants} ensures (\ref{k-fixed-consistent}) in our setting, where
the regular representations of the $G_n$ on the $L^2(M_n)$ are multiplicity 
free.  Thus we have $G_n$--equivariant injections 
\begin{equation}\label{alpha0}
\alpha_{m,n,\lambda^*}: \cH_{n,\lambda^*}^{K_n} \to 
\cH_{m,\lambda^*}^{K_m},\ \ \ \ \
\cH_{n,\lambda}^{K_n} = w_{n,\lambda^*}\C \text{ and } 
\cH_{m,\lambda}^{K_n} = w_{m,\lambda^*}\C
\end{equation}
as in (\ref{scale-invars}), and
\begin{equation}\label{sym-ztild}
\widetilde{\zeta}_{m,n,\lambda}: 
\cH_{n,\lambda}\otimes (w_{n,\lambda^*}\C)
\to \cH_{m,\lambda}\otimes (w_{m,\lambda^*}\C)
\text{ by } f_{u,w_n,n} \mapsto c_{m,n,\lambda}
\left (\tfrac{\deg\pi_{m,\lambda}}{\deg\pi_{n,\lambda}}\right )^{1/2}f_{u,w_m,m}
\end{equation}
as in (\ref{scale-invars1}).  Now as in Section \ref{sec2e},
$\left \{\cH_{n,\lambda} \otimes (w_{n,\lambda^*}\C) , \ \
\widetilde{\zeta}_{m,n,\lambda}\right\}$ is isomorphic to a subsystem
of the system
$\left \{\left (\cH_{n,\lambda} \widehat{\otimes}\cH_{n,\lambda}^*\right ),
\zeta_{m,n,\lambda}\right \}$ of (\ref{n-normalize-inj-cpt}).  As a 
result we have

\begin{equation}\label{rightinv5}
\cH_\lambda \otimes w_{\lambda^*}\C := \varinjlim
\left \{(\cH_{n,\lambda} \otimes w_{n,\lambda^*}\C), \ \
\widetilde{\zeta}_{m,n,\lambda}\right \}
\end{equation}
and they fit together under the direct integral (\ref{rightinv}),
which here is reduced to a discrete direct sum, to give us
$L^2(M)$ where $M = \varinjlim M_n$ as follows.
\begin{equation}\label{rightinv6}
L^2(M) := L^2(G)^K =
\varinjlim \left \{L^2(G_n/K_n), \widetilde{\zeta}_{m,n}\right \}.
\end{equation}

We proceed as in Section \ref{sec2e}, but taking advantage of the fact that 
here the $G_n$ are compact.  Define
\begin{equation}\label{symm-quo-lim-coef}
\begin{aligned}
&\cA(\pi_{n,\lambda})^{K_n} = \{\text{finite linear combinations of the } 
  f_{u,v_n,n}
  \text{ where } u \in \cH_{\pi_{n,\lambda}} \text{ and } 
  v_n \in \cH_{\pi_{n,\lambda}}^{K_n}\}, \\
&\nu_{m,n,\lambda}:\cA(\pi_{n,\lambda})^{K_n} \hookrightarrow 
  \cA(\pi_{m,\lambda})^{K_m} \text{ by }
        f_{u,v_n,n} \mapsto f_{u,v_m,m} \text{ where projection }
        p_{m,n,\lambda}(v_m) = v_n.
\end{aligned}
\end{equation}
Thus Lemma \ref{res-lemma} says:
If $f \in \cA(\pi_{n,\lambda})^{K_n}$ then 
$\nu_{m,n,\lambda}(f)|_{G_n/K_n} = f$.

The ring of regular functions on $M_n = G_n/K_n$ is
$\cA(G_n/K_n) := \cA(G_n)^{K_n}
= \sum_\lambda \cA(\pi_{n,\lambda})$ and the $\nu_{m,n,\lambda}$ sum to
define a direct system $\{\cA(G_n/K_n),\nu_{m,n}\}$.  Its limit is
\begin{equation}\label{symm-quo-lim-reg}
\cA(G/K) := \cA(G)^K =
\varinjlim \{\cA(G_n/K_n),\nu_{m,n}\}.
\end{equation}
As before, the maps of the direct system $\{\cA(G_n/K_n),\nu_{m,n}\}$
are inverse to restriction of functions, so $\cA(G/K)$ is a $G$--submodule
of the inverse limit $\varprojlim \{\cA(G_n/K_n), \text{ restriction}\}$.

Each
$\cA(G_n/K_n)$ is a dense subspace of $L^2(G_n/K_n)$ but, because the 
$\nu_{m,n}$ distort the Hilbert space structure,
$\cA(G/K)$ does not sit naturally as a subspace of $L^2(G/K)$.  Thus we
use the $G_n$--equivariant maps 
\begin{equation}\label{sym-rel-map-sys}
\widetilde{\eta}_{n,\lambda}: \cA(\pi_{n,\lambda})^{K_n} \to
\cH_{\pi_n}\widehat{\otimes}(w_{n,\lambda^*}\C) \text{ by }
f_{u,w_{n,\lambda},n} \mapsto c_{n,1,\lambda}
\sqrt{\deg \pi_{n,\lambda}}\, f_{u,w_{n,\lambda},n}.
\end{equation}
where $c_{m,n,\lambda} = ||p_{m,n,\lambda}(w_{m,\lambda})||$.
Now Proposition \ref{quo-comparison} specializes to
\begin{proposition}\label{sym-quo-comparison}
The maps $\widetilde{\zeta}_{m,n,\lambda}$ of {\rm (\ref{sym-ztild})},
$\nu_{m,n,\lambda}$ of {\rm (\ref{symm-quo-lim-coef})} and
$\widetilde{\eta}_{n,\lambda}$ of {\rm (\ref{sym-rel-map-sys})}
satisfy $$(\widetilde{\eta}_{m,\lambda}\circ \nu_{m,n,\lambda})(f_{u,v,n}) =
(\widetilde{\zeta}_{m,n,\lambda}\circ\widetilde{\eta}_{n,\lambda})(f_{u,v,n})$$
for $f_{u,v,n} \in \cA(\pi_{n,\lambda})^{K_n}$.  Thus they inject
the direct system $\{\cA(G_n)^{K_n}, \nu_{m,n}\}$
into the direct system $\{L^2(G_n)^{K_n},\widetilde{\zeta}_{m,n}\}$.
That map of direct systems defines a $G$--equivariant injection
$$
\widetilde{\eta}: \cA(G/K) \to L^2(G/K)
$$
with dense image.  In particular $\widetilde{\eta}$ defines a pre Hilbert space
structure on $\cA(G/K)$ with completion isometric to $L^2(G/K)$.
\end{proposition}

\section{Limits of Heisenberg Commutative Spaces.}
\label{sec4}
\setcounter{equation}{0}

We now turn to a class of commutative spaces $G_n/K_n$ on which a
nilpotent subgroup of $G_n$ acts transitively.  Then that nilpotent
group must be the nilradical $N_n$ of $G_n$, the isotropy subgroup
$K_n$ must be a compact group of automorphisms of $N_n$, and $G_n$
must be the semidirect product $N_n\rtimes K_n$.  See \cite[Chapter 13]{W3}
for an exposition.  These spaces usually are weakly symmetric, but they
are not symmetric.  Nevertheless they are accessible because $N_n$ is
very similar in structure to the Heisenberg group, and the theory is
modelled on the Heisenberg group case.  Here in Section \ref{sec4} we
study the Heisenberg group case, and we examine more general
commutative nilmanifolds in Section \ref{sec5}.

\subsection{Regular Functions on Heisenberg Groups.}\label{sec4a}
We first consider the case where $G_n$ is the Heisenberg group
$$
H_n = \Im\C + \C^n \text{ with group composition } (z,w)(z',w')
	= (z + z' + \Im \langle w, w' \rangle , w + w').
$$
There are two sorts of irreducible unitary representations.  The ones that
annihilate the center $Z = \Im\C$ of $H_n$ are the ($1$--dimensional)
unitary characters on the vector group $H_n/Z \cong \R^{2n}$.  The ones
that are nontrivial on $Z$, say $\pi_{n,t}$ with central character
$\zeta_t(z,w) = e^{tz}$, are specified by the nonzero real number $t$ and
are realized on the Foch 
space $\cH_{n,t}$ of entire holomorphic functions $f: \C^n \to \C$ such
that $\int_{\C^m} |f(w)|^2 e^{-|t||w|^2} dw < \infty$ where $dw$ is Lebesgue 
measure.  The representation is 
$[\pi_t(z,v)f](w) = e^{tz \pm t\Im(w- v/2)\cdot v} f(w-v)$
where $\pm $ is the sign of $t/|t|$.
\medskip

For each real $t \ne 0$ the representation $\pi_{n,t}$ is square integrable
(modulo the center of $H_n$) and has formal degree $|t|^n$.
\medskip

For each multi-index
$\um = (m_1, \dots , m_n)$, $m_i \geqq 0$,  we have the monomial
$w^{\um} = w_1^{m_1}\dots w_n^{m_n}$.
For $t \ne 0$ there is a constant $c_{n,t} > 0$ such that the 
$\mu_\um := w^\um/\sqrt{\um!} := (w_1^{m_1} \dots w_n^{m_n})/\sqrt{(m_1! \dots m_n!)}$
satisfy 
$$
\int \mu_\um(w)\overline{\mu_{\um'}(w)}\exp(-|t||w|^2) dw =
\begin{cases} 0 \text{ for } \um \ne \um',\\
1/c_{n,t} \text{ for } \um = \um.
\end{cases}
$$
Now we normalize the inner product on $\cH_{\pi_{n,t}} = \cH_{n,t}$ by
\begin{equation}\label{innerfock}
\langle f,f'\rangle = c_{n,t}\int f(w) \overline{f'(w)} \exp(-|t||w|^2) dw.
\end{equation}
Then the monomials  $\mu_\um$ form a complete orthonormal set in 
$\cH_{\pi_{n,t}}$.
\medskip

In Definition \ref{def-l2-dirint} of direct integral now,
take $(Y,\cM,\tau) = (\R\setminus \{0\},\cM,|t|^ndt)$ where $\cM$ consists
of the Borel sets.  For each multiindex $\um$ define
$s_\um : Y: \to \bigcup \cH_{n,t}$ by $s_\um(t) = e^{-|t|}\mu_m$. 
Then $\langle s_\um(t) , s_{\um'}(t)\rangle_{\cH_{n,t}} = 
e^{-2|t|}\delta_{\um,\um'} \in L^1(Y,\tau)$.  In fact 
\begin{equation}\label{gamma}
\int_{-\infty}^\infty \langle s_\um(t) , s_\um(t)\rangle_{\cH_{n,t}} |t|^ndt
= 2\int_0^\infty e^{-2t} t^n dt = 2^{-n}\int_0^\infty e^{-2t} 2^{n+1}t^ndt
= 2^{-n}\Gamma(n+1) = n!/2^n.
\end{equation}
Thus (\ref{dirintmaps}) is satisfied and we have the Hilbert space
of (\ref{dirintspace}):
\begin{equation}\label{left-int}
\cH_n = \int_{-\infty}^\infty \cH_{n,t} |t|^ndt \text{ defined by }
(\R\setminus \{0\},\cM,|t|^ndt) \text{ and the } s_\um(t) = e^{-|t|}\mu_\um.
\end{equation}
Here $|t|^ndt$ is the Plancherel measure for the Heisenberg group $H_n$.
Since $|t|^ndt$ and $dt$ are mutually absolutely continuous on
$\R\setminus \{0\}$ we can reformulate (\ref{left-int}) as
$\cH_n = \int_{-\infty}^\infty \cH_{n,t} \, dt$, which is independent of $n$,
for purposes of taking direct limits.
\medskip 

For multiindices $\um$ and $\um'$ define
$r_{\um,\um'} : Y: \to \bigcup \cH_{n,t}\widehat{\otimes} \cH_{n,t}^*$ by 
$r_{\um,\um'}(t) = e^{-|t|}\mu_\um \otimes \overline{\mu_{\um'}}$.  The 
$\mu_\um \otimes \overline{\mu_{\um'}}$
form a complete orthonormal set in $\cH_{n,t}\widehat{\otimes} \cH_{n,t}^*$
so $\langle r_{\um,\um'}(t) , r_{\um'',\um'''}(t)\rangle_{\cH_{n,t}} =
e^{-2|t|}\delta_{\um,\um''}\delta_{\um',\um'''}$, which is in $L^1(Y,\tau)$ as
before because
$
\int_{-\infty}^\infty \langle r_{\um,\um'}(t) , r_{\um,\um'}(t)\rangle_{\cH_{n,t}} |t|^ndt
= 2\int_0^\infty e^{-2t} t^n dt 
= 2^{-n}\Gamma(n+1) =n!/2^n.
$
Since the $1$--dimensional representations in $\widehat{H_n}$
form a set of Plancherel measure zero, and $|t|^ndt$ is the restriction of
Plancherel measure to $\{\pi_{n,t} \mid t \ne 0\}$ this gives us
\begin{equation}\label{full-int}
L^2(H_n) = \int_{-\infty}^\infty (\cH_{n,t}\widehat{\otimes} \cH_{n,t}^*)|t|^ndt
\end{equation}
where the direct integral is
defined by $(\R\setminus \{0\},\cM,|t|^ndt) \text{ and the }
r_{\um,\um'}(t) = e^{-|t|}\mu_\um \otimes \overline{\mu_{\um'}}$.
Now we carry this over to rational functions.

\begin{definition}\label{reg-fn-heis-0}
For each real polynomial $p(t)$ define 
$p_{\um,\um'}(t) = e^{-|t|}p(t)\mu_\um \otimes \overline{\mu_{\um'}}$.
Let $\cA(H_n)$ denote the set of all finite linear combinations of the
$p_{\um,\um'}$ with multiplication 
$(p'_{\um,\um'}\times p''_{\uell,\uell'})(t) = [p'(t)p''(t)]
(\mu_{\um+\uell}\otimes \overline{\mu_{\um'+\uell'}})$.
The elements of $\cA(H_n)$ are the {\rm regular functions}
on $H_n$ and $\cA(H_n)$ is the {\rm ring of regular functions} on $H_n$.
Note that $\cA(H_n)$ is isomorphic to the ring of polynomials on
the complex extension $\C + \C^n$ of the underlying vector space structure
of $H_n = \Im\C + \C^n$.
\end{definition}

\begin{lemma}\label{reg-fn-heis-1}
The ring $\cA(H_n)$ of regular functions on $H_n$
is a dense subspace of $L^2(H_n)$.
\end{lemma}
\noindent {\bf Proof.}  The pointwise inner product
$\langle p_{\um,\um'}(t), r_{\um'',\um'''}(t) \rangle_{\cH_{n,t}}
= e^{-2|t|}p(t)\delta_{\um,\um''}\delta_{\um',\um'''}$.  Computing as before,
$\int_{-\infty}^\infty e^{-2|t|} |t^k| |t|^ndt = (k+n)!/2^{k+n}$ for
integral $k \geqq 0$.  If $p(t) = \sum_0^d p_kt^k$ this shows
$$
||\langle p_{\um,\um'}(t), r_{\um,\um'}(t) \rangle_{\cH_{n,t}}||_{L^1(\R\setminus \{0\},\cM,|t|^ndt)}
= \int_{-\infty}^\infty |e^{-2|t|} p(t)| |t|^ndt \leqq
\sum_0^d |p_k|(k+n)!/2^{k+n} < \infty,
$$ 
so each $\langle p_{\um,\um'}(t), r_{\um'',\um'''}(t) \rangle_{\cH_{n,t}} \in
L^1(\R\setminus \{0\},\cM,|t|^ndt)$.  
This proves $\cA(H_n) \subset L^2(H_n)$.
\medskip

Let $f \in L^2(H_n)$ orthogonal to $\cA(H_n)$.  Denote $f_t(w) = f(t,w)$.  
and $r_{\um,\um',t} = r_{\um,\um'}(t)$.
By (\ref{full-int}) and Fubini, $f_t$ is orthogonal to  
every $r_{\um,\um',t}$ a.e. $t \in \R$.
For fixed $t$, the $r_{\um,\um',t}$ form a complete orthogonal set in
$\cH_{n,t}\widehat{\otimes} \cH_{n,t}^*$.  Thus $f_t = 0$ a.e. $t \in \R$.
Now $f = 0$ in $L^2(H_n)$.  We have proved that the subspace $\cA(H_n)$
is dense in $L^2(H_n)$.
\hfill $\square$

\subsection{Functions on the Infinite Heisenberg Group.}
\label{sec4b}
Passage to the limit is easy for Heisenberg groups.  The group
inclusions are given by $(z,w) \mapsto (z,w)$, identity on the center 
$\Im \C$ and the usual $\C^n \to \C^m$ by $(w_1, \dots , w_n) \mapsto
(w_1, \dots , w_n, 0, \dots , 0)$ on its complement.  Also, and this
is a key point, the Hilbert space 
$\cH_{n,t}$ sits naturally in $\cH_{m,t}$ as the closed
span of the monomials $\mu_\um := w^\um/\sqrt{\um!}$ for which the exponents
$m_{n+1} = \dots = m_m = 0$.  That gives us the equivariant isometric
inclusions $\zeta_{m,n}: \cH_{n,t} \to \cH_{m,t}$, sending $\mu_\um$ (as an
element of $\cH_{n,t}$) to $\mu_\um$ (as an element of $\cH_{m,t}$).  Now we
have
\begin{equation}\label{lim-heis}
\begin{aligned}
&H_\infty = \Im\C + \C^\infty = \varinjlim H_n: 
	\text{ infinite Heisenberg group,}\\
&\cH_t = \varinjlim \cH_{n,t}: \text{ Hilbert space with complete 
	orthonormal set } \{\mu_\um = w^\um/\sqrt{\um!}\},\\
&\pi_t = \varinjlim \pi_{n,t}: \text{ irreducible unitary representation of }
	H_\infty \text{ on } \cH_t.
\end{aligned}
\end{equation}

The space $\cE_{n,t}$ of matrix coefficients of $\cH_{n,t}$ is spanned
by the functions 
$f_{\uell,\um;n,t}: g \mapsto \langle \mu_\uell, \pi_{n,t}(g)\mu_\um\rangle$, 
These coefficients belong to the Hilbert space
$$
L^2(H_n/\Im\C;e^t) := \{f:H_n \to \C \mid |f| \in L^2(H_n/\Im\C)
        \text{ and } f(z,w) = e^{-tz}f(0,w)\}
$$
with inner product $\langle f,f'\rangle = \int_{\C^n} f(z,w)\overline{f'(z,w)}
dw$.
\medskip

Since $\pi_{n,t}$ has formal degree $|t|^n$ the orthogonality
relations say that the inner product in $\cE_{n,t}$ is 
$\langle f_{\uell,\um;n,t}, f_{\uell',\um';n,t}\rangle = |t|^{-n}$ if $\uell = \uell'$
and $\um = \um'$, $0$ otherwise.  Now the $|t|^{n/2}f_{\uell,\um;n,t}$
form a complete orthonormal set in $\cE_{n,t}$, and $\cE_{n,t}$ consists of
the functions $\Phi_{n,t,\varphi}$ given by
\begin{equation}
\Phi_{n,t,\varphi}(h) =
\sum_{\uell,\um} \varphi_{\uell,\um}(t) |t|^{n/2} f_{\uell,\um;n,t}(h)
\end{equation}
where $\sum_{\uell,\um} |\varphi_{\uell,\um}(t)|^2 < \infty$.
Thus $L^2(H_n)$ consists of all functions 
\begin{equation}
\Psi_{n,\varphi}(h): = \int_{-\infty}^\infty \Phi_{n,t,\varphi}(h)|t|^n\, dt
= \int_{-\infty}^\infty \left ( {\sum}_{\uell,\um}
\varphi_{\uell,\um}(t) |t|^{n/2}f_{\uell,\um;n,t}(h) \right ) |t|^n\, dt
\end{equation}
such that the $\varphi_{\uell,\um}: \R \to \C$ are measurable,
$\sum_{\uell,\um} |\varphi_{\uell,\um}(t)|^2 < \infty$ a.e. $t$, and
$\sum_{\uell,\um} |\varphi_{\uell,\um}(t)|^2 \in L^1(\R,|t|^ndt)$.
Note that
$||\Psi_{n,\varphi}||^2_{L^2(H_n)} 
= {\sum}_{\uell,\um} ||\varphi_{\uell,\um}||^2_{L^2(\R,|t|^{n/2}dt)}$.
\medskip

When $m \geqq n$ we view an $n$--tuple $\um$ as an $m$--tuple by
appending $m-n$ zeroes; then $w^{\um}$ has the same meaning
as function on $\C^n$ and on $\C^m$.  Thus the coefficient function
$f_{\uell,\um;n,t}: H_n \to \C$ is the restriction of 
$f_{\uell,\um;n,t}: H_m \to \C$.   Now, as in (\ref{incl-gen}) and
(\ref{normalize-inj-cpt}), we have isometric $(H_n\times H_n)$--equivariant 
injections
\begin{equation}\label{zeta-nilp}
\zeta_{m,n,t}: \cE_{n,t} \to \cE_{m,t} \text{ by } 
\zeta_{m,n,t}(|t|^{n/2}f_{\uell,\um;n,t}) = |t|^{m/2}f_{\uell,\um;n,t} \text{ and }
\zeta_{m,n}(\Psi_{n,\varphi}) = \Psi_{m,|t|^{(n-m)/2}\varphi}
\end{equation}
of $\cE_{n,t}$ into $\cE_{m,t}$ and $L^2(H_n)$ into $L^2(H_m)$.  
Specifically, $\zeta_{m,n,t}$ maps a complete
orthonormal set in $\cE_{n,t}$ to an orthonormal set in $\cE_{m,t}$
and $\zeta_{m,n}$ passes to the direct integral.  As expected, the result
is the multiplicity--free left/right regular representation of 
$H_\infty \times H_\infty$ on
\begin{equation}\label{lim-zeta-nilp}
L^2(H_\infty) := \varinjlim \{L^2(H_n), \zeta_{m,n}\}.
\end{equation}

Now look back to Definition \ref{reg-fn-heis-0}, and define
\begin{equation}\label{reg-inf-heis}
\cA(H_\infty) = \varinjlim \cA(H_n) = \{\text{finite linear combinations
of the } p_{\um,\um'}(t) = e^{-|t|}p(t)w^\um \overline{w^{\um'}}\}.
\end{equation}
Thus $\cA(H_\infty)$ consists of the finite linear combinations of the
$e^{-|t|} t^k \mu_\um \, \overline{\mu_{\um'}}$, and (for $n$
sufficiently large so that $\mu_\um,\, \mu_{\um'} \in \cH_{n,t}$)
we have
$||e^{-|t|} t^k \mu_\um \, \overline{\mu_{\um'}}||^2_{L^2(H_n)}
= (k+n)!/2^{k+n}$.  Consider the maps of (\ref{map-sys}):
$$
\eta_{n,t}: \cA(\pi_{n,t}) \to \cH_{n,t}\widehat{\otimes}\cH_{n,t}^* 
\text{ given by } 
\eta_{n,t}(e^{-|t|} p(t) \mu_\um \, \overline{\mu_{\um'}}) =
e^{((\tfrac{n}{2}-1)|t|)} p(t) \mu_\um \, \overline{\mu_{\um'}}.
$$
Then Proposition \ref{comparison} tells us that 
\begin{proposition}\label{heis-comparison}
The maps $\eta_{m,t}$ satisfy 
$\eta_{m,t}\circ\zeta_{m,n,t} = \zeta_{m,n,t}\circ\eta_{n,t}$
on $\cA(\pi_{n,t})$ and send the direct system $\{\cA(H_n)\}$
into the direct system $\{L^2(H_n)\}$.  That direct system map defines
an $(H_\infty \times H_\infty)$--equivariant injection
$\eta:\cA(H_\infty) \to L^2(H_\infty)$ with dense image.  In particular
$\eta$ defines a pre Hilbert space structure on $\cA(H_\infty)$ with completion
isometric to $L^2(H_\infty)$.
\end{proposition}

\subsection{Heisenberg Nilmanifolds.}\label{4c}
We first carry the results of Section \ref{sec2e} and
\ref{sec4b} over to commutative spaces
$G_n/K_n$ where $G_n$ is the semidirect product $H_n\rtimes K_n$  of a 
Heisenberg group with a compact group of automorphisms.
The concrete results in this section will require that $K_n$ be
connected and that its action on $\C^n$ be irreducible.
\medskip

The classification goes as follows for the cases where $K_n$ is 
connected and is irreducible 
on $\C^n$.  Carcano's Theorem (\cite{Ca}; or see 
\cite[Theorem 4.6]{BJR} or \cite[Theorem 13.2.2]{W3}) says
that $(G_n,K_n)$ is a Gelfand pair if and only if
the representation of $(K_n)_{_\C}$, on polynomials on $\C^n$, is
multiplicity free.  Those groups were classified by Ka\v c 
\cite[Theorem 3]{K} in another context.
His list (as formulated in \cite[(13.2.5)]{W3}) is

{\small
\begin{equation} \label{kac-table}
\begin{tabular}{|r|l|l|l|l|}\hline
\multicolumn{5}{| c |}{Irreducible connected groups $K_n \subset U(n)$
multiplicity free on polynomials on $\C^n$} \\
\hline \hline
& Group $K_n$ & Group $(K_n)_{_\C}$ & Acting on & Conditions on $n$ \\ \hline
1 & $SU(n)$ & $SL(n;\C)$ & $\C^n$ & $n \geqq 2$ \\ \hline
2 & $U(n)$ & $GL(n;\C)$ & $\C^n$ & $n \geqq 1$ \\ \hline
3 & $Sp(m)$ & $Sp(m;\C)$ & $\C^n$ & $n = 2m$ \\ \hline
4 & $U(1) \times Sp(m)$ & $\C^* \times Sp(m;\C)$ & $\C^n$ & $n = 2m$ \\ \hline
5 & $U(1) \times SO(n)$ & $\C^* \times SO(n;\C)$ & $\C^n$ & $n \geqq 2$
        \\ \hline
6 & $U(m)$ & $GL(m;\C)$ & $S^2(\C^m)$ & $m \geqq 2, \
        n = \tfrac{1}{2}m(m+1)$  \\ \hline
7 & $SU(m)$ & $SL(m;\C)$ & $\Lambda^2(\C^m)$ & $m$ odd,
        $n = \tfrac{1}{2}m(m-1)$  \\ \hline
8 & $U(m)$ & $GL(m;\C)$ & $\Lambda^2(\C^m)$ &
        $n = \tfrac{1}{2}m(m-1)$  \\ \hline
9 & $SU(\ell) \times SU(m)$ & $SL(\ell;\C) \times SL(m;\C)$ &
        $\C^\ell \otimes \C^m$ & $n = \ell m, \ \ell \ne m$ \\ \hline
10 & $U(\ell) \times SU(m)$ & $GL(\ell;\C) \times SL(m;\C)$ &
        $\C^\ell \otimes \C^m$ & $n = \ell m$ \\ \hline
11 & $U(2) \times Sp(m)$ & $GL(2;\C) \times Sp(m;\C)$ &
        $\C^2 \otimes \C^{2m}$ & $n = 4m$ \\ \hline
12 & $SU(3) \times Sp(m)$ & $SL(3;\C) \times Sp(m;\C)$ &
        $\C^3 \otimes \C^{2m}$ & $n = 6m$ \\ \hline
13 & $U(3) \times Sp(m)$ & $GL(3;\C) \times Sp(m;\C)$ &
        $\C^3 \otimes \C^{2m}$ & $n = 6m$ \\ \hline
14 & $U(4) \times Sp(4)$ & $GL(4;\C) \times Sp(4;\C)$ &
        $\C^4 \otimes \C^8$ & $n = 32$ \\ \hline
15 & $SU(m) \times Sp(4)$ & $SL(m;\C) \times Sp(4;\C)$ &
        $\C^m \otimes \C^8$ & $n = 8m, \ m \geqq 3$ \\ \hline
16 & $U(m) \times Sp(4)$ & $GL(m;\C) \times Sp(4;\C)$ &
        $\C^m \otimes \C^8$ & $n = 8m, \ m \geqq 3$ \\ \hline
17 & $U(1) \times Spin(7)$ & $\C^* \times Spin(7;\C)$ & $\C^8$ &
        $n = 8$ \\ \hline
18 & $U(1) \times Spin(9)$ & $\C^* \times Spin(9;\C)$ & $\C^{16}$ &
        $n = 16$ \\ \hline
19 & $Spin(10)$ & $Spin(10;\C)$ & $\C^{16}$ & $n = 16$ \\ \hline
20 & $U(1) \times Spin(10)$ & $\C^* \times Spin(10;\C)$ &
        $\C^{16}$ & $n = 16$ \\ \hline
21 & $U(1) \times G_2$ & $\C^* \times G_{2,\C}$ & $\C^7$ & $n=7$ \\ \hline
22 & $U(1) \times E_6$ & $\C^* \times E_{6,\C}$ & $\C^{27}$ & $n=27$ \\ \hline
\end{tabular}
\end{equation}
}
Now we have the direct systems 
{\normalsize
\begin{equation} 
\label{jaw-table}
\begin{tabular}{|c|l|l|l|}\hline
\multicolumn{4}{| c |}{Direct systems $\{(H_n\rtimes K_n,K_n)\}$ of Gelfand pairs,}
\\
\multicolumn{4}{| c |}{$K_n$ connected and irreducible on $\C^n$} \\
\hline \hline
 & Group $K_n$ & Acting on & Conditions on $n$ \\ \hline
1 & $SU(n)$ & $\C^n$ & $n \geqq 2$ \\ \hline
2 & $U(n)$ & $\C^n$ & $n \geqq 1$ \\ \hline
3 & $Sp(m)$ & $\C^n$ & $n = 2m$ \\ \hline
4 & $U(1) \times Sp(m)$ & $\C^n$ & $n = 2m$ \\ \hline
5a & $U(1) \times SO(2m)$ & $\C^{2m}$ & $n= 2m \geqq 2$
        \\ \hline
5b & $U(1) \times SO(2m+1)$ & $\C^{2m+1}$ & $n = 2m+1\geqq 3$
        \\ \hline
6 & $U(m)$ & $S^2(\C^m)$ & $m \geqq 2, \
        n = \tfrac{1}{2}m(m+1)$  \\ \hline
7 & $SU(m)$ & $\Lambda^2(\C^m)$ & $m$ odd,
        $n = \tfrac{1}{2}m(m-1)$  \\ \hline
8 & $U(m)$ & $\Lambda^2(\C^m)$ &
        $n = \tfrac{1}{2}m(m-1)$  \\ \hline
9 & $SU(\ell) \times SU(m)$ &
        $\C^\ell \otimes \C^m$ & $n = \ell m, \ \ell \ne m$ \\ \hline
10 & $S(U(\ell) \times U(m))$ &
        $\C^\ell \otimes \C^m$ & $n = \ell m$ \\ \hline
11 & $U(2) \times Sp(m)$ &
        $\C^2 \otimes \C^{2m}$ & $n = 4m$ \\ \hline
12 & $SU(3) \times Sp(m)$ &
        $\C^3 \otimes \C^{2m}$ & $n = 6m$ \\ \hline
13 & $U(3) \times Sp(m)$ &
        $\C^3 \otimes \C^{2m}$ & $n = 6m$ \\ \hline
15 & $SU(m) \times Sp(4)$ &
        $\C^m \otimes \C^8$ & $n = 8m, \ m \geqq 3$ \\ \hline
16 & $U(m) \times Sp(4)$ &
        $\C^m \otimes \C^8$ & $n = 8m, \ m \geqq 3$ \\ \hline
\end{tabular}
\end{equation}
}
\hskip -.15 cm
In each case the direct system $\{K_n\}$ is both strict and parabolic.
(We separated entry 5 of Table \ref{kac-table} into entries 5a and 5b of 
Table \ref{jaw-table} in order to have the parabolic property.)
\medskip

We now suppose that $\{K_n\}$ is one of the  strict parabolic direct
system given by the rows of 
Table \ref{jaw-table}.
\medskip

As $U(n)$ acts on $H_n = \Im\C + \C^n$ by $k: (z,v) \mapsto (z,kv)$ it
preserves the equivalence class of each of the square integrable representations
$\pi_{n,t}$ of $H_n$.  The Mackey obstruction vanishes and $\pi_{n,t}$
extends to a unitary representation $\widetilde{\pi_{n,t}}$ of 
$H_n\rtimes U(n)$ on $\cH_{n,t}$.  See \cite[Section 4]{W1} for a geometric
proof.  We will also write $\widetilde{\pi_{n,t}}$
for its restriction, the extension of $\pi_{n,t}$ to a unitary representation
of $G_n = H_n\rtimes K_n$.  
\medskip

Let $\kappa_{n,\lambda}$ denote the irreducible representation of $K_n$
of highest weight $\lambda$, and $\widetilde{\kappa_{n,\lambda}}$
its extension to $G_n = H_n\rtimes K_n$ which annihilates $H_n$.
The corresponding representation space $\cF_{n,\lambda}$ is a finite
dimensional Hilbert space.  Denote $\pi_{n,t,\lambda} = 
\widetilde{\pi_{n,t}}\otimes \widetilde{\kappa_{n,\lambda}}$.  Then
$\cH_{n,t,\lambda} := \cH_{n,t}\otimes \cF_{n,\lambda}$ 
is its representation space. 
Fix an orthonormal basis $\{u_i\}$ of $\cF_{n,\lambda}$.  Then 
$\{\mu_\um \otimes u_i\}$ is a complete orthonormal set in 
$\cH_{n,t,\lambda}$, and we have matrix coefficients 
$$
f_{\uell,\um,i,j;n,t,\lambda}(h,k) = \langle (\mu_\uell \otimes u_i), 
((\widetilde{\pi_{n,t}}\otimes 
\widetilde{\kappa_{n,\lambda}})(h,k))(\mu_\um\otimes u_j)\rangle.
$$
The formal degree $\deg \pi_{n,t,\lambda} = |t|^n\deg(\kappa_{n,\lambda})$,
so the $|t|^{n/2}\deg(\kappa_{n,\lambda})^{1/2} f_{\uell,\um,i,j;n,t,\lambda}$ 
form a complete orthonormal set in the space $\cE_{n,t,\lambda}
= \cH_{n,t,\lambda} \widehat{\otimes} \cH_{n,t,\lambda}^*$ of matrix 
coefficient functions.  As in (\ref{incl-gen}), (\ref{normalize-inj-cpt}) and
(\ref{zeta-nilp}) we have isometric $(G_n\times G_n)$--equivariant
injections
{\small
\begin{equation}\label{k-zeta-nilp1}
\zeta_{m,n,t,\lambda}: \cE_{n,t,\lambda} \to \cE_{m,t,\lambda} \text{ by }
\zeta_{m,n,t,\lambda}((|t|^n\deg\kappa_{n,\lambda})^{1/2}
	f_{\uell,\um,i,j;n,t,\lambda}) = 
(|t|^m\deg\kappa_{m,\lambda})^{1/2} f_{\uell,\um,i,j;m,t,\lambda}
\end{equation}
}
Integrate with respect to $t$ and sum on $\lambda$ to construct
isometric $(G_n\times G_n)$--equivariant injections
$\zeta_{m,n}: L^2(G_n) \to L^2(G_m)$.

\begin{theorem}\label{lim-sd}
For $n > 0$ let $K_n$ be a compact connected subgroup of $\Aut(H_n)$
such that $\{K_n\}$ is a strict parabolic direct system.
Define $G_n = H_n\rtimes K_n$, $G = \varinjlim \{G_n\}$ and
$K = \varinjlim \{K_n\}$.  Note $G = H_\infty \rtimes K$.
Then $\{L^2(G_n),\zeta_{m,n}\}$ is a strict direct system of Hilbert spaces
in which the  maps $\zeta_{m,n}: L^2(G_n) \to L^2(G_m)$ are
$(G_n\times G_n)$--equivariant unitary injections.  
That gives us the left/right regular representation of $G\times G$
on the Hilbert space $L^2(G): = \varinjlim \{L^2(G_n), \zeta_{m,n}\}$.
Further, that left/right regular representation 
is the multiplicity--free 
$\int_{-\infty}^\infty (\pi_{t,\lambda}\boxtimes \pi_{t,\lambda}^*)\,\,dt$ where
$\pi_{t,\lambda} := \varinjlim \pi_{n,t,\lambda}$.
\end{theorem}

Now the construction of $\cA(G)$ follows the lines of (\ref{reg-inf-heis}),
with properties relative to $L^2(G)$ as in Proposition \ref{heis-comparison}.
As in Definition \ref{reg-fn-heis-0} we define
\begin{equation}\label{reg-fn-heis-2}
\cA(G) = \varinjlim \cA(G_n) \text{ where }
\cA(G_n) = \{\text{finite linear combinations of the }
	e^{-|t|}p(t)f_{\uell,\um,i,j;n,t,\lambda}\}
\end{equation}
where $p(t)$ is a real polynomial in $t$.  The argument of 
Lemma \ref{reg-fn-heis-1} shows that $\cA(G_n)$ is a dense subspace
of $L^2(G_n)$.  The maps of (\ref{map-sys}) in this setting are
$$
\eta_{n,t,\lambda}: \cA(\pi_{n,t,\lambda}) \to 
\cH_{n,t,\lambda}\widehat{\otimes}\cH_{n,t,\lambda}^*
\text{ by }
\eta_{n,t,\lambda}(e^{-|t|} p(t) f_{\uell,\um,i,j;n,t,\lambda}) =
e^{(\tfrac{n}{2}-1)|t|} \sqrt{\deg \kappa_{n,\lambda}}
\,\, p(t) f_{\uell,\um,i,j;n,t,\lambda}.
$$
Then Proposition \ref{comparison} tells us that
\begin{proposition}\label{heis-comparison-again}
The maps $\eta_{m,t,\lambda}$ satisfy $\eta_{m,t,\lambda}\circ
\zeta_{m,n,t,\lambda} = \zeta_{m,n,t,\lambda}\circ\eta_{n,t,\lambda}$
on $\cA(\pi_{n,t,\lambda})$ and send the direct system $\{\cA(G_n)\}$
into the direct system $\{L^2(G_n)\}$.  That system map defines
an $\left ((H_\infty\rtimes K) \times (H_\infty\rtimes K)\right )$--equivariant 
injection $\eta:\cA(H_\infty\rtimes K) \to L^2(H_\infty\rtimes K)$ with 
dense image.  In particular $\eta$ defines a pre Hilbert space structure on 
$\cA(H_\infty\rtimes K)$ with completion
isometric to $L^2(H_\infty\rtimes K)$.
\end{proposition}

Recall our working hypothesis that $\{K_n\}$ is one of the $16$ systems of
Table \ref{jaw-table}.  
Since $(G_n, K_n)$ is a Gelfand pair with $K_n$ irreducible
on $\C^n$, Carcano's Theorem \cite{Ca} says that the action of $K_n$ on
the polynomial ring $\C[C^n]$ is multiplicity free, and it picks out the
right $K_n$--invariants in $L^2(G_n)$, as follows.

\begin{lemma}\label{kntriv}
Recall the notation of {\rm Section \ref{sec3a}} for $\Xi_n$, $\lambda$ and
$\kappa_{n,\lambda} \in \widehat{K_n}$.
Define $\widetilde{\kappa_{n,\lambda}} \in \widehat{G_n}$ by 
$\widetilde{\kappa_{n,\lambda}}(hk) = \kappa_{n,\lambda}(k)$ for
$h \in H_n$ and $k \in K_n$.  Then $\pi_{n,t,\lambda} := 
\widetilde{\pi_{n,t}}\otimes \widetilde{\kappa_{n,\lambda}}$ has a
nonzero $K_n$--fixed vector if and only if $\kappa_{n,\lambda}^*$ occurs as a
subrepresentation of $\widetilde{\pi_{n,t}}|_{K_n}$, and in that
case the space of $K_n$--fixed vectors has dimension $1$.
\end{lemma}

\noindent {\bf Proof.} This is essentially the argument in 
\cite[Section 14.5A]{W1}.  Decompose $\widetilde{\pi_{n,t}}|_{K_n} = 
\sum_{\gamma \in \widetilde{K_n}} m_\gamma\ \gamma$.  Carcano's Theorem
(\cite{Ca}, or see \cite[Theorem 13.2.2]{W3}) says that each $m_\gamma$ is
either $0$ or $1$.  The 
$K_n$--fixed vectors of $\widetilde{\kappa} \otimes \widetilde{\pi_{n,t}}$ all
occur in $\kappa \otimes (m_{\kappa^*}\kappa^*)$, and they form a space
of dimension $m_{\kappa^*}$.  The assertion follows. 
\phantom{XXXXXXXXXX} \hfill $\square$
\medskip

We view $L^2(G_n/K_n)$ as the space
of right--$K_n$--invariant functions in $L^2(G_n)$.  With
Lemma \ref{kntriv} in mind we set
$$
\widehat{K_n}^\dagger = \{\kappa_{n,\lambda} \in \widehat{K_n} \mid
\kappa_{n,\lambda}^* \text{ occurs in the space of polynomials on } \C^n\} .
$$
Recall that $\cF_{n,\lambda}$ denotes the representation space of 
$\kappa_{n,\lambda}$.
Given $\kappa_{n,\lambda} \in \widehat{K_n}^\dagger$ the right $K_n$--invariant
in $\C[\C^n]\otimes \cF_{n,\lambda}^*$ is 
$\int_{K_n}\sum_i (b_i\otimes \kappa_{n,\lambda}^*(b_i^*))dk$ 
where $\{b_i\}$ is a basis of the
$\kappa_{n,\lambda}$--subspace of $\C[\C^n]$ and $\{b_i^*\}$ is the dual
basis of 
$\cF_{n,\lambda}^*$.  Normalize it to a unit vector $w_{n,t,\lambda}$.
Then the (left regular) representation of $G_n$ on $L^2(G_n/K_n)$ is 
equivalent to the representation
$\sum_{\kappa_{n,\lambda} \in \widehat{K_n}^\dagger} \int_{-\infty}^\infty
\widetilde{\pi_{n,t}}\otimes \widetilde{\kappa_{n,\lambda}}\, dt$ of
$G_n$ on
$
\sum_{\kappa_{n,\lambda} \in \widehat{K_n}^\dagger} \int_{-\infty}^\infty
(\cH_{n,t,\lambda}\otimes w_{n,t,\lambda}\C)\, dt.
$

\begin{proposition}\label{res-inv}
If $m \geqq n$ and $\kappa_{n,\lambda} \in \widehat{K_n}^\dagger$
then $\kappa_{m,\lambda} \in \widehat{K_m}^\dagger$.  In that case
inclusion $\C[\C^n] \hookrightarrow \C[\C^m]$ of polynomials maps
the highest weight $\lambda$ space for $\kappa_{n,\lambda}$ 
onto the the highest weight $\lambda$ space for $\kappa_{m,\lambda}$.
\end{proposition}

\noindent {\bf Proof.}
The group $K_n$ acts on $\C^n$ by some representation
$\gamma_n$, so the representation of $K_n$ on polynomials of
degree $d$ is the symmetric power $S^d(\gamma_n^*)$.  Thus we can compute
the set $X_{n,d}$ of highest weights of $K_n$ on the space $P_{n,d}$ of 
polynomials of degree $d$ on $\C^n$.  Running through the $16$ cases of 
Table \ref{jaw-table} we see that $X_{n,d} \subset X_{m,d}$.  
For example (Line 3 of Table \ref{jaw-table}) the representation of $Sp(m)$ on
polynomials of degree $q$ in $C^{2m}$ is the irreducible representation
with highest weight $q\xi_{m,1}$, and (Lines 5a and 5b of 
Table \ref{jaw-table})
the representation of $U(1)\times SO(n)$ on polynomials of degree $q$
in $C^{n}$ is the tensor product of the $-q^{th}$ power of the usual 
representation of $U(1)$ by scalars on $\C^{n}$ with 
the multiplicity--free sum of irreducible representations of $SO(n)$ of
highest weights $\{\xi_{n,1}, 3\xi_{n,1}, 5\xi_{n,1}, ... , q\xi_{n,1}\}$ 
if $q$ is odd,
$\{0\xi_{n,1}, 2\xi_{n,1}, 4\xi_{n,1}, ... , q\xi_{n,1}\}$ if $q$ is even.
\medskip

Now let $\lambda \in X_{n,d}$.  Let $v_{n,\lambda}$ denote a (nonzero)
highest weight $\lambda$ vector for $\gk_n$ in $P_{n,d}$, and
similarly let $v_{m,\lambda}$ denote a (nonzero) highest weight $\lambda$ 
vector for $\gk_m$.  Divide up the variables of $\C^m$ to 
$\{w_1, \dots , w_n\}$ for $\C^n$ and $\{z_{n+1}, \dots , z_m\}$ for its 
complement in $\C^m$.  Express $v_{m,\lambda} = \sum_{A,B}b_{A,B}w^Az^B$
where each term has total degree $|A|+|B| = d$.  Note that $K_n$ treats
the $z_i$ as constants.  Evaluating the $z_i$ at arbitrary constant values
$C = (c_{n+1}, \dots , c_m)$ we have a highest weight $\lambda$ vector for
$\gk_n$. By Carcano's Theorem it is a multiple of $v_{n,\lambda}$. In other 
words $v_{m,\lambda}|_{\{z = C\}} = m_{_C}v_{n,\lambda}$.  The terms
$b_{A,B}\uw^A\uz^B$ with $\uz$--degree $|B| > 0$ yield evaluations of 
$\uw$--degree $|A| < d$, and cannot contribute to any $m_{_C}v_{n,\lambda}$.
Now $b_{A,B}\uw^A\uz^B = 0$ whenever $|B| > 0$.  This shows that
$v_{m,\lambda}$ is a homogeneous polynomial of degree $d$ in the $w_j$,
as is $v_{n,\lambda}$.  We conclude that $v_{m,\lambda}$ is a nonzero
multiple of $v_{n,\lambda}$.
\hfill $\square$

\begin{corollary}\label{invar-nest}
Condition {\rm (\ref{k-fixed-consistent})} is satisfied for the
direct systems $\{K_n\}$ of {\rm Table \ref{jaw-table}}.
\end{corollary}

\noindent {\bf Proof.} Retain the notation $X_{n,d}$ for those
$\lambda$ such that $\kappa_{n,\lambda}$ occurs on the space $P_{n,d}$
of polynomials of degree $d$ on $\C^n$.  If $\lambda \notin X_{n,d}$
there are no nonzero $K_n$--invariant vectors in $\cH_{n,t,\lambda}$,
so the assertion is vacuous.  Now assume $\lambda \in X_{n,d}$
and choose an orthonormal basis $\{x_1, \dots , x_{q(n)}\}$ of the 
representation space for $\kappa_{n,\lambda}$ in $P_{n,d}$.  According 
to Proposition \ref{res-inv} that representation space is contained in 
the representation space for $\kappa_{m,\lambda}$ in $P_{m,d}$, so the 
latter has an orthonormal basis
$\{x_1, \dots , x_{q(n)}, x_{q(n)+1}, \dots , x_{q(m)}\}$.  Let
$\{x^*_1, \dots , x^*_{q(n)}, x^*_{q(n)+1}, \dots , x^*_{q(m)}\}$ and
$\{x^*_1, \dots , x^*_{q(n)}\}$ 
be the corresponding dual bases of $\cF_{m,\lambda}$ and $\cF_{n,\lambda}$.
The $K_n$--invariant vectors in $\cH_{n,t,\lambda}$ are the multiples
of $\sum_1^{q(n)} x_i\otimes x_i^*$, and the $K_m$--invariant vectors in
$\cH_{m,t,\lambda}$ are the multiples of $\sum_1^{q(m)} x_i\otimes x_i^*$.
The adjoint of unitary inclusion is orthogonal projection, which sends
$\sum_1^{q(m)} x_i\otimes x_i^*$ to $\sum_1^{q(n)} x_i\otimes x_i^*$.
\phantom{XXXXXXXXX}\hfill $\square$
\medskip

Let $w_{n,\lambda}$ be a right--$K_n$--fixed unit vector
in the highest weight $\lambda$ subspace of $L^2(H_n\rtimes K_n)$, for each
$\kappa_{n,\lambda} \in \widehat{K_n}^\dagger$.  Note that $w_{n,\lambda}$
does not depend on $t$.  Proposition \ref{res-inv}
says that the inclusion $\cE_{n,t,\lambda} \hookrightarrow \cE_{m,t,\lambda}$
maps $w_{n,\lambda}$ to a nonzero multiple of $w_{m,\lambda}$.  Given
$\lambda$ we recursively choose the $w_{n,\lambda}$ so that 
\begin{equation}\label{def-c}
w_{m,\lambda} = c_{m,n,t,\lambda}w_{n,\lambda} + x \text{ with }
x \perp \cE_{n,t,\lambda} \text{ with } 0 < c_{n,t,\lambda} \leqq 1.
\end{equation}
Note $\cH_{n,\lambda^*}^{K_n} = w_{n,\lambda^*}\C$ and
$\cH_{m,\lambda^*}^{K_n} = w_{m,\lambda^*}\C$.
Now we have $(H_n\rtimes K_n)$--equivariant injections
\begin{equation}\label{alpha1}
\alpha_{m,n,t,\lambda^*}: \cH_{n,t,\lambda^*}^{K_n} \to
\cH_{m,t,\lambda^*}^{K_m},\ \ \ \
\alpha_{m,n,t,\lambda^*}(w_{n,\lambda^*}) = c_{m,n,t,\lambda}w_{m,\lambda^*},
\end{equation}
as in (\ref{scale-invars}), and
\begin{equation}\label{heis-rel-scale}
\widetilde{\zeta}_{m,n,t,\lambda}:
\cH_{n,t,\lambda}\otimes (w_{n,t,\lambda^*}\C)
\to \cH_{m,t,\lambda}\otimes (w_{m,t,\lambda^*}\C)
\text{ defined by } f \mapsto c_{m,n,t,\lambda}\zeta_{m,n,t,\lambda}(f)
\end{equation}
as in (\ref{quo-lim-coef}) and (\ref{quo-lim-coef1}).  Now as in Section 
\ref{sec2e},
$\left \{\cH_{n,t,\lambda} \otimes (w_{n,\lambda^*}\C) , \ \
\widetilde{\zeta}_{m,n,t,\lambda}\right\}$ is isomorphic to a subsystem
of the system
$\left \{\left (\cH_{n,t,\lambda} \widehat{\otimes}\cH_{n,t,\lambda}^*\right ),
\zeta_{m,n,t,\lambda}\right \}$ of (\ref{n-normalize-inj-cpt}).  As a
result we have
\begin{equation}\label{rightinv9}
\cH_{t,\lambda} \otimes w_{\lambda^*}\C := \varinjlim
\left \{\cH_{n,t,\lambda}\otimes w_{n,\lambda^*}\C, \ \
\widetilde{\zeta}_{m,n,t,\lambda}\right \}
\end{equation}
and they fit together under the direct integral (\ref{rightinv})
to give us $L^2((H_n\rtimes K_n)/K_n)$ as follows.
\begin{equation}\label{rightinvz}
L^2((H_\infty\rtimes K)/K) := L^2(H_\infty\rtimes K)^K =
\varinjlim \left \{L^2((H_n\rtimes K_n)/K_n), 
\widetilde{\zeta}_{m,n}\right \}.
\end{equation}
Combining Theorem \ref{lim-sd}, Lemma \ref{kntriv} and Corollary 
\ref{invar-nest} we have

\begin{theorem}\label{heis-case}
Let $\{G_n,K_n)\}$ be one of the direct systems of
{\rm Table \ref{jaw-table}}. Define $G_n = H_n\rtimes K_n$,
$G = \varinjlim G_n$ and $K = \varinjlim K_n$.  Then {\rm (\ref{rightinvz})} is
a unitary direct system whose limit
Hilbert space is $G$--isometric to
$L^2(G)^K$, and the natural unitary representation of $G$ on 
$L^2(G/K) = L^2(G)^K$ is multiplicity free.
\end{theorem}

Now we turn to regular functions.  As in (\ref{reg-inf-heis}) we define
\begin{equation}\label{reg-fin-heis-k}
\begin{aligned}
&\cA(H_n\rtimes K_n) := \{\text{finite linear combinations of the }
        e^{-|t|}f_{\uell,\um,i,j;n,t,\lambda} \text{ in }
        \cE_{n,t,\lambda}\}, \\
& \cA((H_n\rtimes K_n)/K_n) := \cA(H_n\rtimes K_n)^{K_n} =
	\cA(H_n\rtimes K_n)\cap \cE_{n,t,\lambda}^{K_n}, \\
&\nu_{m,n,t,\lambda}:\cA(\pi_{n,t,\lambda})^{K_n} \hookrightarrow
  \cA(\pi_{m,t,\lambda})^{K_m} \text{ by }
        f_{u,v_n,n} \mapsto f_{u,v_m,m} \text{ where }
        p_{m,n,t,\lambda}(v_m) = v_n
\end{aligned}
\end{equation}
Now we have direct systems and their limits
\begin{equation}\label{reg-inf-heis-k}
\begin{aligned}
& \cA(H_\infty\rtimes K) = \varinjlim \{\cA(H_n\rtimes K_n), \zeta_{m,n}\}
\text{ where } \zeta_{m,n,t,\lambda} : \cE_{n,t,\lambda} 
\hookrightarrow \cE_{m,t,\lambda} \text{  (\ref{k-zeta-nilp1}), and}\\
&\cA((H_\infty\rtimes K)/K) = 
\varinjlim \{\cA((H_n\rtimes K_n)/K_n), \nu_{m,n}\}
\text{ where } \nu_{m,n,t,\lambda} : \cE_{n,t,\lambda}^{K_n}
\hookrightarrow \cE_{m,t,\lambda}^{K_m} \text{  (\ref{reg-fin-heis-k})}.
\end{aligned}
\end{equation}

As before, each $\cA((H_n\rtimes K_n)/K_n)$ 
is a dense subspace of $L^2((H_n\rtimes K_n)/K_n)$.
In order to pass this comparison to the limit we use the maps
\begin{equation}\label{k-heis-rel-map-sys}
\widetilde{\eta}_{n,t,\lambda}: \cA(\pi_{n,t,\lambda})^{K_n} \to
\cH_{\pi_{n,t,\lambda}}\otimes(w_{n,t,\lambda^*}\C) \text{ by }
f \mapsto c_{n,1,\lambda}\, |t|^{n/2}\, \sqrt{\deg \kappa_{n,\lambda}}\, f.
\end{equation}

\begin{proposition}\label{k-heis-quo-comparison}
The $\widetilde{\eta}_{n,t,\lambda}$
satisfy $(\widetilde{\eta}_{m,t,\lambda}\circ
\nu_{m,n,t,\lambda})(f) =
(\widetilde{\zeta}_{m,n,\lambda}\circ\widetilde{\eta}_{n,t,\lambda})(f)$
for $f \in \cA(\pi_{n,t,\lambda})^{K_n}$  Thus they inject
the direct system $\{\cA((H_n\rtimes K_n)/K_n), \nu_{m,n}\}$
of regular functions into the direct system 
$\{L^2((H_n\rtimes K_n)/K_n),\widetilde{\zeta}_{m,n}\}$ of square integrable
functions.  That map of direct systems defines an 
$(H_\infty\rtimes K)$--equivariant injection
$$
\widetilde{\eta}: \cA((H_\infty\rtimes K)/K) \to L^2((H_\infty\rtimes K)/K)
$$
with dense image.  In particular $\widetilde{\eta}$ defines a pre Hilbert space
structure on $\cA((H_\infty\rtimes K)/K)$ with completion isometric 
to $L^2((H_\infty\rtimes K)/K)$.
\end{proposition}

\section{Extension to Commutative Nilmanifolds.}\label{sec5}
\setcounter{equation}{0}
The results of Section \ref{sec4} depend on four basic facts.  First, the
$\pi_{n,t}$ are determined by their central character.  Second,
we have good models $\cH_{n,t}$ for the representation spaces, such that
$n$ does not appear explicitly in the formulae for the actions of the group
elements.  Third, the injections $H_n \hookrightarrow H_m$ restrict
to isomorphisms $Z_n \cong Z_m$ of the centers.  And fourth,
we have complete information on the Plancherel measure for the $H_n$.
In this section we consider a somewhat larger class of nilpotent
direct systems that satisfy these conditions.
In this section we extend our Heisenberg group considerations to nilpotent 
Lie groups with square integrable representations, 
following the general lines of \cite{W3} and \cite{W5}.  

\subsection{Square Integrable Nilpotent Groups.}\label{sec5a}
Here is a quick summary of harmonic analysis 
for connected simply connected groups that admit square integrable
representations.  See \cite{MW} for details, 
\cite[Section 14.2]{W3} for an exposition.
Let $N$ be a connected, simply connected nilpotent Lie group and $\gn$ its Lie
algebra.  Decompose $\gn = \gz + \gv$ and $N = Z\exp(\gv)$ where $\gz$
is the center of $\gn$.  Then $Z = \exp(\gz)$ is the center of $N$.  We
say that an irreducible unitary representation $\pi$ of $N$ is
{\sl square integrable} if its coefficient functions $f_{u,v}(g) =
\langle u, \pi(g)v\rangle$ satisfy $|f_{u,v}| \in L^2(N/Z)$.  In that case
$\pi$ is
determined by its central character, $\pi = \pi_t$ where $t \in \gz^*$
and the central character is $\exp(\zeta) \mapsto e^{it(\zeta)}$.  In
terms of geometric quantization, $\pi_t$ corresponds to the coadjoint orbit
in $\gn^*$ consisting of all linear functionals on $\gn$ whose restriction
to $\gz$ is $t$.  Further, the
antisymmetric bilinear form $b_t(\xi,\eta) = t([\xi,\eta])$ on $\gv$
is nondegenerate, and (up to a positive constant that depends only on
normalizations of Haar measures) the formal degree of $\pi_t$
is $|\Pf(b_t)|$, where $\Pf(b_t)$ is the Pfaffian\footnote{Strictly speaking,
$\Pf(b_t)$ depends on a choice of basis of $\gv$, for a basis change
of determinant $a_t$ multiplies $\det b_t|_{\gv \times \gv}$ by
$\det a_t^2$ and multiplies $\Pf(b_t)$ by $\det a_t$.} of
$b_t: \gv \times \gv \to \R$.  In fact, if $\pi_s$ is the representation
of $N$ that corresponds to $\Ad^*(N)s \subset \gn^*$, then
$\pi_s$ is square integrable if and only if $\Pf(b_s) \ne 0$.
In any case, 
\begin{equation}\label{pf-poly}
\Pf(b_t) \text{ is a polynomial function of } t,\,\, \Pf(b_t) = P(t)
\text{, and } P(0) = 0.
\end{equation} 
Again, up
to a constant that depends on normalizations, $|\Pf(b_t)|$ is the
Plancherel density.  It follows that if one irreducible unitary
representation of $N$ is square integrable then Plancherel--almost--all are.
In the case of the Heisenberg group $H_n$, where we identified $\gz^*$
with $\R$, the Pfaffian corresponding to $\pi_{n,t}$ is $t^n$.
\medskip

When we are dealing with a sequence $\{N_n\}$ if square integrable
nilpotent groups, we have to keep track of the polynomials (\ref{pf-poly}),
so we will write
\begin{equation}\label{pf-poly-n}
P(n,t) = \Pf(b_{n,t}): \text{ corresponding to the group } N_n\, .
\end{equation}

The point of this, from the viewpoint of commutative spaces, is that many
Gelfand pairs are of the form $(N\rtimes K,K)$ where $N$ is a connected
simply connected Lie group, $K$ is a compact subgroup of $\Aut(N)$, and
$N$ has square integrable representations.  See \cite[Theorem 14.4.3]{W3}.
This is simplified by the $2$--step Nilpotent Theorem
\cite[Theorem 13.1.1]{W3} of Benson-Jenkins-Ratcliff and Vinberg, which says
that $N$ must be abelian or $2$--step nilpotent.  In a certain
sense representations treat those groups as Heisenberg groups:

\begin{lemma} \label{is-heis1} {\rm (\cite[Lemma 14.4.1]{W3})}
Let $N$ be a connected simply connected $2$--step nilpotent
Lie group with $1$--dimensional center.  Then $N$ is isomorphic to the
Heisenberg group $H_n$ where $n = \tfrac{1}{2}(\dim_{_\R} \gn - 1)$, and in
particular $N$ has square integrable representations.
\end{lemma}

\begin{proposition} \label{is-heis2} {\rm (\cite[Proposition 14.4.2]{W3})}
Let $N$ be a connected simply connected $2$--step nilpotent
Lie group.  Let $f \in \gn^*$ such that $f|_\gz \ne 0$.
Denote $\gw_f = \{z \in \gz \mid f(z) = 0\}$ and $W_f := \exp(\gw_f)$.
Then

{\rm 1.} $W_f$ is a closed subgroup of $Z$, hence a closed normal subgroup
of $N$.

{\rm 2.} The functional $f$ is the pullback of a linear functional
$f' \in (\gn/\gw_f)^*$ and is nonzero on the central subalgebra $\gz/\gw_f$
of $\gn/\gw_f$.

{\rm 3.} The representation $[\pi_f]$ is the pullback to $N$ of the
class $[\pi_{f'}] \in \widehat{N/W_f}$.

{\rm 4.} If the representation $[\pi_f]$ is square integrable then
$[\pi_{f'}]$ is square integrable, and in that case
$N/W_f$ has center $Z/W_f$ and is isomorphic to a Heisenberg group
$H_n$ where $n = \tfrac{1}{2}\dim(\gn/\gz)$.
\end{proposition}

We now consider a strict direct system $\{N_n\}$ of $2$--step nilpotent
connected, simply connected Lie groups that have square integrable
representations, where the inclusions $\gn_n \to \gn_m$ map
the center $\gz_n \hookrightarrow \gz_m$ and the complement
$\gv_n \hookrightarrow \gv_m$ in decompositions $\gn_n = \gz_n + \gv_n$.
Then the direct limit algebra $\gn := \varinjlim \gn_n$ has center
$\gz := \varinjlim \gz_n$ and $\gn = \gz + \gv$ where
$\gv = \varinjlim \gv_n$.  On the group
level, $Z = \varinjlim Z_n = \exp(\gz)$ is the center of
$N := \varinjlim N_n$ and we have $N = Z\exp(\gv)$.
\medskip

We further assume that the dimensions $\dim Z_n$ of the centers are
bounded.  Since they are non--decreasing we may assume that they are
eventually constant.  Passing to a cofinal sequence,
\begin{equation} \label{samecenter}
\gn_n \hookrightarrow \gn_m \text{ maps } \gz_n \cong \gz_m.
\end{equation}
Under that identification we write $\gz$ for all the $\gz_n$,\,\,
$\gz^*$ for all the $\gz_n^*$,\,\, and $Z$ for all the $Z_n$.
\medskip

Let $t \in \gz^*$.  Write $b_{n,t}$ for the
bilinear form  $(\xi,\eta) \mapsto t([\xi,\eta])$ on
$\gv_n$.  Then $t$ corresponds to a square integrable representation 
$\pi_{n,t}$ of $N_n$ just when the Pfaffian $\Pf(b_{n,t}) \ne 0$.
For purposes of comparing the Pfaffians as $n$ varies, we note that
$\Pf(b_{n,t})$ is specified by $t$ and a basis of $\gv_n$, so we
simply assume that these bases are nested in the sense that the basis
of $\gv_{n+1}$ consists the basis of $\gv_n$ together with some elements
that are $b_{n+1,t}$--orthogonal to $\gv_n$.  Thus, if
$\Pf(b_{n+1,t}) \ne 0$ then $\Pf(b_{n,t}) \ne 0$.  The converse fails
in general, but the following lemma deals with the possibility that
$\Pf(b_{n,t}) \ne 0 = \Pf(b_{m,t})$.  It depends on the fact \cite{MW}
that each $\Pf(b_{n,t})$ is a polynomial function on $\gz^*$.

\begin{lemma}\label{v-basis}
Let $\ga_n \in \gz^*$ denote the zero set of $\Pf(b_{n,t})$ and
set $\ga = \bigcup \ga_n$.  Then $\ga$ is a set of Lebesgue measure
zero in $\gz^*$.
\end{lemma}
\noindent {\bf Proof.}  Since $N_n$ has square integrable representations,
the Pfaffian $\Pf(b_{n,t})$ is a nontrivial polynomial function of
$t \in \gz^*$, so $\ga_n$ is a finite union of lower--dimensional
subvarieties of $\gz_n^*$.  Now the set $\ga$ is a countable union of
sets $\ga_n$ of Lebesgue measure zero.  \hfill $\square$
\medskip

Define $T = \{t \in \gz^* \mid \text{ each } \Pf(b_{n,t}) \ne 0\}$. 
So $\gz^* \setminus \ga$.  For every $t \in T$ and every index $n$ we have
a square integrable representation $\pi_{n,t} \in \widehat{N_n}$.
Let $t \in T$, $\gw_t = \{z \in \gz \mid t(z) = 0\}$ and $W_t = \exp(\gw_t)$.
Then $W_t$ is closed in $Z$, $N_n/W_t$
is isomorphic to a Heisenberg group $H_{d(n)}$, and
$\pi_{n,t}$ factors
through to the square integrable representation of $N_n/W_t$ with
central character $e^{it}$.  As $t$ varies in $T$ the $\pi_{n,t}$ act on the same Fock space $\cH_{d(n),t}$,
$d(n) = \tfrac{1}{2}\dim \gv_n$, by formulae independent of $d(n)$.
\medskip

We normalize the inner products on the $\cH_{d(n),t}$ as before,
so the $\mu_\um$ form a complete orthonormal set, and realize the
space $\cE_{n,t} = \cH_{d(n),t}\widehat{\otimes} \cH_{d(n),t}^*$ of
matrix coefficients as the closed span of the functions $f_{\uell,\um;n,t}
: g \mapsto \langle \mu_\uell,\pi_{n,t}(g)\mu_\um \rangle$,
as in Section \ref{sec4}.  The orthogonality relations say that the
inner product on $\cE_{n,t}$ is given by
$\langle f_{\uell,\um;n,t}, f_{\uell',\um';n,t}\rangle
= |\Pf(b_{n,t})|^{-1}$ if $\uell = \uell'$ and
$\um = \um'$, and is $0$ otherwise.  Now the 
$|\Pf(b_{n,t})|^{1/2}f_{\uell,\um;n,t}$
form a complete orthonormal set in $\cE_{n,t}$, and as before
$\cE_{n,t}$ consists of the functions $\Phi_{n,t,\varphi}$ on $H_{d(n)}$
given by
\begin{equation}
\Phi_{n,t,\varphi}(h) =
\sum_{\uell,\um} \varphi_{\uell,\um}(t) |\Pf(b_{n,t})|^{1/2}f_{\uell,\um;n,t}(h)
\text{ with } \sum_{\uell,\um} |\varphi_{\uell,\um}(t)|^2 < \infty.
\end{equation}
Now
$L^2(N_n) = \int_{\gz_n^*} \cE_{n,t} |\Pf(b_{n,t})|
dt = \int_T \cE_{n,t} |\Pf(b_{n,t})| dt$.  It consists of all functions
$\Psi_{n,\varphi}$ defined by
\begin{equation}\label{c1}
\Psi_{n,\varphi}(h) =  \int_{\gz_n^*} \Phi_{n,t,\varphi}(h)
|\Pf(b_{n,t})| dt 
= \int_T \left ( {\sum}_{\uell,\um}\varphi_{\uell,\um}(t)
        |\Pf(b_{n,t})|^{1/2} f_{\uell,\um;n,t} \right ) |\Pf(b_{n,t})| dt
\end{equation}
such that the functions $\varphi_{\uell,\um}: \gz_n^* \to \C$
are measurable with $\sum_{\uell,\um}
|\varphi_{\uell,\um}(t)|^2 < \infty$ for almost all $t \in T$ and
$\sum_{\uell,\um} |\varphi_{\uell,\um}(t)|^2 \in L^1(\gz_n^*, |\Pf(b_{n,t})| dt)$.
The norms are
\begin{equation}\label{c2}
\begin{aligned}
||\Psi_{n,\varphi}||^2_{L^2(N_n)}
        &= \int_T ||\Phi_{n,t,\varphi}||^2_{\cE_{n,t}}
                |\Pf(b_{n,t})| dt 
        = \int_T \left ( {\sum}_{\uell,\um}
                |\varphi_{\uell,\um}(t)|^2 \right ) |\Pf(b_{n,t})| dt \\
        &= {\sum}_{\uell,\um}
         ||\varphi_{\uell,\um}||^2_{L^2(\gz^*,|\Pf(b_{n,t})|dt)}\ \ .
\end{aligned}
\end{equation}
The left/right representation of
$N_n\times N_n$ on $\cE_{n,t}$ is the exterior tensor product
$\pi_{n,t} \boxtimes \pi_{n,t}^*$; it is irreducible and the left/right
regular representation of $N_n\times N_n$ on $L^2(N_n)$ is
the multiplicity free unitary representation
$\int_{\gz^*} (\pi_{n,t} \boxtimes \pi_{n,t}^*) |\Pf(b_{n,t})|dt$.
\medskip

Let $m \geqq n$.  Then
$|\Pf(b_{n,t})|^{1/2} f_{\uell,\um;n,t} \mapsto |\Pf(b_{m,t})|^{1/2} f_{\uell,\um;m,t}$ 
defines an equivariant isometric injection
$\Phi_{n,\varphi}(t) \mapsto \Phi_{m,\varphi}(t)$ of
$\cE_{n,t}$ into $\cE_{m,t}$.  The norm computation just above
gives
\begin{equation}\label{c3}
||\Psi_{m,|\Pf(b_{n,t})/\Pf(b_{m,t})|^{1/2}\,\varphi}||^2_{L^2(N_m)}
= ||\Psi_{n,\varphi}||^2_{L^2(N_n)}.
\end{equation}
Thus we have an $(N_n\times N_n)$--equivariant isometric injection
\begin{equation}\label{c4}
\zeta_{m,n}: L^2(N_n) \to L^2(N_m) \text{ by }
        \zeta_{m,n}(\Psi_{n,\varphi}) =
          \Psi_{m,|\Pf(b_{n,t})/\Pf(b_{m,t})|^{1/2}\varphi}.
\end{equation}
On the level of coefficients it is given by
$\zeta_{m,n}(\Phi_{n,t,\varphi}) =
\Phi_{m,t,|\Pf(b_{n,t})/\Pf(b_{m,t})|^{1/2}\varphi}$.  In other words
$\zeta_{m,n}$ sends the function $\sum_{\uell,\um} \varphi_{\uell,\um}(t)
                |\Pf(b_{n,t})|^{1/2}f_{\uell,\um;n,t}$
on $N_n$ to the function on $N_m$ given by
$$
\sum_{\uell,\um} (|\Pf(b_{n,t})/\Pf(b_{m,t})|^{1/2}
        \varphi_{\uell,\um}(t))(|\Pf(b_{m,t})|^{1/2}f_{\uell,\um;m,t}) 
= \sum_{\ul,\um} \varphi_{\uell,\um}(t) |\Pf(b_{n,t})|^{1/2}f_{\uell,\um;t}.
$$
The result is
\begin{theorem} \label{nilp-inj}
There is a strict direct system $\{L^2(N_n), \zeta_{m,n}\}$ of $L^2$ spaces.
The direct system maps $\zeta_{m,n}: L^2(N_n) \to L^2(N_m)$ are
$(N_n\times N_n)$--equivariant unitary injections.  Let $\Pi_n$ denote the
left/right regular representation of $N_n\times N_n$ on $L^2(N_n)$
and let $N = \varinjlim N_n$.  Then we have a well defined Hilbert
space $L^2(N): = \varinjlim \{L^2(N_n), \zeta_{m,n}\}$ and a 
multiplicity--free unitary
representation $\Pi = \varinjlim \Pi_n$ of $N \times N$ on
$L^2(N)$.  
\end{theorem}

\subsection{Semidirect Product Groups.}\label{sec5b}
We need some structural results from \cite[Section 7]{W5} for a 
strict direct system $\{K_n, \varphi_{m,n}\}$ of compact connected Lie
groups and a consistent family $\{\gamma_n\}$ of representations of the
$K_n$ on a fixed finite dimensional vector space $\gz$.  We'll use 
that information to extend Theorems \ref{lim-sd}
and \ref{heis-case} to a larger family of strict direct systems of
nilmanifold Gelfand pairs.
\medskip

As just indicated, $\{K_n, \varphi_{m,n}\}$ is a strict direct system of
compact connected Lie groups.  Denote $K = \varinjlim \{K_n, \varphi_{m,n}\}$
and let $\{\gamma_n\}$ be a consistent family of unitary representations of the
$K_n$ on a fixed finite dimensional real vector space $\gz$.  The 
$U_n = \gamma_n(K_n)$ form an increasing sequence of compact
connected subgroups of dimension $\leqq (\dim\gz)^2$ in the unitary
group $U(\gz)$, so from some index on they are all the same compact connected
subgroup $U$ of $U(\gz)$.  Now we truncate the index set and
assume that the $\gamma_n$ all have the
same image $U$ in the orthogonal group of $\gz$.
\medskip

Let $K_n^\dagger$ denote the
identity component of the kernel of $\gamma_n$.  Then
$\varphi_{m,n}(K_n^\dagger) \subset K_m^\dagger$, so we have
$K^\dagger = \varinjlim \{K_n^\dagger, \varphi_{m,n}|_{K_n^\dagger}\}$,
and $K^\dagger$ is the identity component of the kernel of $\gamma$.
Since $K_n$ is compact and connected, and $K_n^\dagger$ is a closed
connected normal subgroup, $K_n$ has another closed connected normal
subgroup $L_n$ such that $K_n$ is locally isomorphic to the direct
product $K^\dagger_n \times L_n$.

\begin{proposition}\label{choose-center}
One can choose the groups $L_n$ so that $\varphi_{m,n}(L_n) = L_m$
for $m \geqq n \ggg 0$.
\end{proposition}

Let $L = \varinjlim\{L_n,\varphi_{m,n}|_{L_n}\}$.
We further truncate the index set so that $L_n = \varphi_n^{-1}(L)$ for all 
indices $n$.

\begin{corollary} \label{parab-cor}
Let $t \in \gz$, and let $K_{n,t}$ be its stabilizer in $K_n$.
If one of the direct systems
$\{K_n, \varphi_{m,n}\}$, $\{K_{n,t}, \varphi_{m,n}|_{K_{n,t}}\}$, or
$\{K_n^\dagger, \varphi_{m,n}|_{K_n^\dagger}\}$ is parabolic, then
the other two also are parabolic.
\end{corollary}

\begin{corollary} \label{parab-cor2}
Let $t \in \gz$, and let $K_{n,t}$ be its stabilizer in $K_n$.
Suppose that the direct system $\{K_n, \varphi_{m,n}\}$ is
parabolic.  Then there are natural isometric injections
$\cF_{n,t,\lambda} \hookrightarrow \cF_{m,t,\lambda}$ for $m\geqq n$,
from the highest weight $\lambda$ representation space of $K_{n,t}$
to that of $K_{m,t}$, and corresponding $L^2$--isometric injections
$f \mapsto \bigl ( (\deg \kappa_{m,t,\lambda})/(\deg\kappa_{n,t,\lambda}) 
\bigr )^{1/2}f$ on spaces of coefficient functions.
\end{corollary}

\begin{corollary}\label{compact-on-z}
Let $t \in \gz$, and let $K_{n,t}$ be its stabilizer in $K_n$.
Then $L := \varinjlim L_n$ is compact, $K = K^\dagger L$ where
$K^\dagger := \varinjlim K_n^\dagger$, and
$K$ is locally isomorphic to $K^\dagger \times L$.  In
particular $K$ acts on $\gz$ as a compact linear group and $\gz$ has
a $\gamma(K)$--invariant positive definite inner product.
\end{corollary}

Now we can proceed along the lines of \cite[Section 8]{W5}.  Fix 
a strict direct system $\{(G_n,K_n)\}$ of Gelfand pairs that satisfies
\begin{equation}\label{nil-conditions}
\begin{aligned}
\text{(i) }&G_n = N_n\rtimes K_n \text{, semidirect product, where } N_n
        \text{ is a connected, }
        \text{ simply connected,}\\
	&\text{ nilpotent Lie group with square integrable
         representations and } K_n \text{ is connected,}\\
\text{(ii) }&\text{the } K_n \text{ form a parabolic strict direct system},  \\
\text{(iii) }&\text{the inclusions } \gn_n \hookrightarrow \gn_{n+1} 
	\text{ map centers }
        \gz_n \cong \gz_{n+1} \text{ and map 
         complements } \gv_n \hookrightarrow \gv_{n+1}, \\
\text{(iv) }&\text{for each $n$ the complement $\gv_n$ is
        $\Ad(K_n)$--invariant.}
\end{aligned}
\end{equation}
Using (\ref{nil-conditions}(iii)) we identify each $\gz_n$ with $\gz :=
\varinjlim \gz_n$.  Write $K_n^\dagger$ for the identity component of the 
kernel of the action of $K_n$ on $\gz$.  Since $\Ad(K_n)|_\gz$ is a compact
connected connected group of linear transformations of $\gz$, its dimension
is bounded, so the $\Ad(K_n)|_\gz$ stabilize and we may assume that each
$\Ad(K_n)|_\gz = U$ for some compact connected group $U$ of linear 
transformations of $\gz$.
Proposition \ref{choose-center} gives us complementary closed connected normal
subgroups $L_n \subset K_n$ that inject isomorphically under 
$K_n \hookrightarrow
K_{n+1}$, so each $L_n$ is equal to $L := \varinjlim L_n$.
Thus we have decompositions $K_n = K_n^\dagger \cdot L$ and
$K = K^\dagger\cdot L$ where $K_n^\dagger$ is the kernel of the adjoint
action of $K_n$ on $\gz$, $K = \varinjlim K_n$, and
$K^\dagger = \varinjlim K_n^\dagger$.  For each $n$,
$\Ad_{K_n}$ maps $L = L_n$ onto $U$ with finite kernel.
\medskip

If $t \in \gz^*$ denote $\cO_t = \Ad^*(L)(t)$, and let 
$G_{n,t}$, $K_{n,t}$ and $L_{n,t}$ denote the respective stabilizers of $t$ 
in $G_n$, $K_n$ and $L_n$.  Since
$\Ad^*(G_n)$ acts on $\cO_t$ as the compact group $L$ there is an
invariant measure $\nu_t$ derived from euclidean measure on $\gz^*$.
Given $t \in \gz^*$ its
stabilizers $G_t = \{g \in G \mid \Ad^*(g)t = t\}$, $K_t = K\cap G_t$
and $L_t = L\cap G_t$.  Their pullbacks in $G_n$ are $G_{n,t}$,
$K_{n,t}$ and $L_{n,t}$.   Note that $K_t = K^\dagger \cdot L_t$
and $K_{n,t} = K_n^\dagger \cdot L_{n,t}$.
\medskip

Recall $T = \{t \in \gz^* \mid \text{ each } \Pf(b_{n,t}) \ne 0\}$ and
fix $t \in T$.  Essentially as in the Heisenberg group case, the
square integrable representation $\pi_{n,t}$ extends to a unitary
representation $\widetilde{\pi_{n,t}}$ of $G_{n,t} :=
N_n\rtimes K_{n,t}$ on the same representation space $\cH_{n,t}$.  If
$\kappa_{n,t,\lambda} \in \widehat{K_{n,t}}$ has representation space
$\cF_{n,t,\lambda}$ we write $\widetilde{\kappa_{n,t,\lambda}}$ for its
extension to a representation of $G_{n,t}$ on
$\cF_{n,t,\lambda}$ that annihilates $N_n$.  Then we have the irreducible
unitary representations
\begin{equation}\label{diamond}
\pi_{n,t,\lambda}^\diamondsuit :=
\widetilde{\pi_{n,t}} \otimes \widetilde{\kappa_{n,t,\lambda}} \in
\widehat{G_{n,t}} \text{ with representation space }
\cH_{n,t,\lambda}^\diamondsuit :=
\cH_{n,t} \otimes \cF_{n,t,\lambda}.
\end{equation}
That gives us the unitarily induced representations
\begin{equation}\label{diamond-induced}
\begin{aligned}
\pi_{n,t,\lambda} = &\Ind_{G_{n,t}}^{G_n} (\pi_{n,t,\lambda}^\diamondsuit)
\in \widehat{G_n}\\& \text{ with representation space }
\cH_{n,t,\lambda} := \int_{\cO_t} (\cH_{n,\Ad^*(k)t}
\otimes \cF_{n,\Ad^*(k)t,\Ad^*(k)\lambda})\, d\nu_t(k(t)).
\end{aligned}
\end{equation}
According to the Mackey little group theory, (i) $\pi_{n,t,\lambda}$ is
irreducible, (ii) $\pi_{n,t,\lambda}$ is equivalent to
$\pi_{n,t',\lambda'}$ if and only if $t' \in \cO_t$, say
$t' = \Ad^*(\ell)t$ where $\ell \in L$, and $\Ad^*(\ell)$ carries $\lambda$
to $\lambda'$, and (iii)
Plancherel--almost--all irreducible unitary representations of
$G_n$ are of the form $\pi_{n,t,\lambda}$ where $t \in T$
and $\kappa_{n,t,\lambda} \in \widehat{K_{n,t}}$.
So the left/right regular representation of $G_n$ 
\begin{equation}\label{sum-diamond-induced}
\begin{aligned}
\Pi_n = &\sum_\lambda \int_{\Ad^*(K_n)\backslash \gz^*}
	(\pi_{n,t,\lambda} \boxtimes \pi_{n,t,\lambda}^*) d(\Ad^*(K_n)(t))
\\ &\text{ with representation space } L^2(G_n) =
\sum_\lambda \int_{\Ad^*(K_n)\backslash \gz^*}
   (\cH_{n,t,\lambda} \widehat{\otimes} \cH_{n,t,\lambda}^*) d(\Ad^*(K_n)(t)).
\end{aligned}
\end{equation}
Since $\pi_{n,t,\lambda}^\diamondsuit$ is square integrable 
and $\pi_{n,t,\lambda}$ is irreducible, one knows (\cite[Theorem A.1]{W5})
that $\pi_{n,t,\lambda}$ has a well defined formal degree.
Making use of Corollaries \ref{parab-cor}, \ref{parab-cor2} and
\ref{compact-on-z}, 
we have $G_n$--equivariant isometric injections
\begin{equation}\label{lim-diamond-induced}
\begin{aligned}
\zeta_{m,n}: L^2(G_n) \to L^2(G_m) &\text{ based on the }
\zeta_{m,n,t,\lambda}: \cH_{n,t,\lambda} \widehat{\otimes} \cH_{n,t,\lambda}^*
\to \cH_{m,t,\lambda} \widehat{\otimes} \cH_{m,t,\lambda}^*\\
&\text{ where } \zeta_{m,n,t,\lambda}((\deg \pi_{n,t,\lambda})^{1/2}f) = 
(\deg \pi_{m,t,\lambda})^{1/2}f.
\end{aligned}
\end{equation}
With $G := \varinjlim G_n$ we now have the 
\begin{equation}\label{lim-diamond-induced-1}
\text{left/right regular representation }
\Pi := \varinjlim \Pi_n \text{ of } G \times G \text{ on } 
L^2(G):= \varinjlim \{L^2(G_n),\zeta_{m,n}\}.
\end{equation}

Now we turn to regular functions.  As one might expect,
$\cE_{n,t,\lambda}^\diamondsuit$ means
$\cH_{n,t,\lambda}^\diamondsuit \widehat{\otimes} 
(\cH_{n,t,\lambda}^\diamondsuit)^*$ viewed as
matrix coefficients of $\pi_{n,t,\lambda}^\diamondsuit$.
As in the (\ref{reg-fn-heis-2}) we define
\begin{equation}\label{reg-fin-nil-kk}
\begin{aligned}
&\cA(G_{n,t}) := \{\text{finite linear combinations of the }
e^{-||t||}f_{\uell,\um,i,j;n,t,\lambda} \text{ in }
\cE_{n,t,\lambda}^\diamondsuit\}, \\
& \cA(G_{n,t}) \to \cA(G_{m,t}) \text{ by } 
e^{-||t||}f_{\uell,\um,i,j;n,t,\lambda} \mapsto 
e^{-||t||}f_{\uell,\um,i,j;m,t,\lambda},
\text{ and } \cA(G_t) = \varinjlim \cA(G_{n,t}).
\end{aligned}
\end{equation}
The norm $||t||$ on $\gz^*$ is from a $U$--invariant inner product.
\medskip

The representation space $\cH_{n,t,\lambda}$ of
$\pi_{n,t,\lambda} = \Ind_{G_{n,t}}^{G_n} (\pi_{n,t,\lambda}^\diamondsuit)$
is the space of $L^2$ sections of the Hilbert space bundle
$\H_{n,t,\lambda}^\diamondsuit \to \cO_t = K_n/K_{n,t}$.
We express it as the space 
$(L^2(K_n)\otimes \cF_{n,t,\lambda}\otimes \cH_{n,t})^{K_{n,t}}$ of 
$K_{n,t}$--invariants in $L^2(K_n)\otimes \cF_{n,t,\lambda}\otimes \cH_{n,t}$.
Let $\cH_{n,t}^{\text{poly}}$ denote the space of finite linear
combinations of the $\mu_\um$ in the space $\cH_{n,t}$.  
Then the underlying algebraic space is the space
$\cA(\pi_{n,t,\lambda})^{K_{n,t}} := (\cA(K_n)\otimes \cF_{n,t,\lambda}\otimes 
\cH_{n,t}^{\text{poly}})^{K_{n,t}}$ of $K_{n,t}$--invariants in
$\cA(K_n)\otimes \cF_{n,t,\lambda}\otimes \cH_{n,t}^{\text{poly}}$.
With that and (\ref{reg-fin-nil-kk}) in mind we define
{\small
\begin{equation}\label{reg-fin-nil-kkk}
\begin{aligned}
&\cA(G_n) := \left \{\text{finite linear comb of the }
e^{-||t||}p(t)f_{n,t,\lambda} \text{ where }
f_{n,t,\lambda} \in (\cA(K_n)\otimes \cF_{n,t,\lambda}\otimes
\cH_{n,t}^{\text{poly}})^{K_{n,t}}\right \}, \\
& \cA(G_n) \to \cA(G_m) \text{ by }
e^{-||t||}p(t)f_{\uell,\um,i,j;n,t,\lambda} \mapsto
e^{-||t||}p(t)f_{\uell,\um,i,j;m,t,\lambda},
\text{ and } \cA(G) = \varinjlim \cA(G_n).
\end{aligned}
\end{equation}
}
where the $p(t)$ are real polynomials on $\gz^*$.
The maps corresponding to those of (\ref{map-sys}) are the
$$
\eta_{n,t,\lambda}: \cA(\pi_{n,t,\lambda}) \to \cE_{n,t,\lambda}
\text{ by }
\eta_{n,t,\lambda}(e^{-||t||} p(t) f_{\uell,\um,i,j;t,\lambda}) =
e^{-||t||} \sqrt{\deg \pi_{n,t,\lambda}}
\,\, p(t) f_{\uell,\um,i,j;t,\lambda}
$$
Proposition \ref{comparison} now tells us that
\begin{proposition}\label{nil-comparison-again}
The maps $\eta_{m,t,\lambda}$ satisfy $\eta_{m,t,\lambda}\circ
\zeta_{m,n,t,\lambda} = \zeta_{m,n,t,\lambda}\circ\eta_{n,t,\lambda}$
on $\cA(\pi_{n,t,\lambda})$ and send the direct system $\{\cA(G_n)\}$
into the direct system $\{L^2(G_n),\zeta_{m.n}\}$.  That system map defines
an $(G \times G)$--equivariant
injection $\eta:\cA(G) \to L^2(G)$ with
dense image.  In particular $\eta$ defines a pre Hilbert space structure on
$\cA(G)$ with completion isometric to $L^2(G)$.
\end{proposition}

\subsection{Commutative Quotients.}\label{sec5c}
We modify the results of Section \ref{sec5b} to strict direct systems
of commutative spaces that satisfy (\ref{nil-conditions}).  Here the $L^2$
results are recalled from \cite[Section 9]{W5}.
\begin{theorem} \label{iso-gelfand}
Let $t \in T$.  Then $(G_{n,t},K_{n,t})$ is a Gelfand pair.
In particular $K_{n,t}$ is multiplicity free on $\C[\gv_n]$.
\end{theorem}

We have already used Hilbert bundle model for the induced representation
$\pi_{n,t,\lambda} \in \widehat{G_{n,t}}$, given by
$\pi_{n,t,\lambda} = \Ind_{G_{n,t}}^{G_n}(\pi_{n,t,\lambda}^\diamondsuit)$.
The representation space $\cH_{n,t,\lambda}$
of $\pi_{n,t,\lambda}$ consists of all
$L^2(K_n/K_{n,t})$ sections of the homogeneous bundle
$p: \H_{n,t,\lambda}^\diamondsuit \to G_n/G_{n,t} = K_n/K_{n,t}$ whose typical
fiber is the representation space $\cH_{n,t,\lambda}^\diamondsuit$ of
$\pi_{n,t,\lambda}^\diamondsuit$.  Given $k \in K_n$ we write
$k\cdot \cH_{n,t,\lambda}^\diamondsuit$ for the fiber $p^{-1}(kK_{n,t})$.
Let $u \in \cH_{n,t,\lambda}^\diamondsuit$ be a
$\pi_{n,t,\lambda}^\diamondsuit (K_{n,t})$--fixed unit vector.
Then $u$ belongs to the fiber $1\cdot \cH_{n,t,\lambda}^\diamondsuit$, and
$k\cdot u \in k\cdot \cH_{n,t,\lambda}^\diamondsuit$ depends only on the
coset $kK_{n,t}$.  Define a section
\begin{equation}\label{induced-invariant}
\sigma_u : K_n/K_{n,t} \to  \H_{n,t,\lambda}^\diamondsuit \text{ by }
        \sigma_u(kK_{n,t}) = k\cdot u.
\end{equation}
Then $\sigma_u$ is a $\pi_{n,t,\lambda}(K_n)$--invariant unit vector in
the Hilbert space $\cH_{n,t,\lambda}$.  (We will also
write $\varphi_u$ for the corresponding function
$G_n \to \cH_{n,t,\lambda}^\diamondsuit$ such that
$\varphi_u(gg_t) = \pi_{n,t,\lambda}^\diamondsuit(g_t)^{-1}(\varphi_u(g))$
for $g \in G_n$ and $g_t \in G_{n,t}$.)
Conversely if $\sigma$ is a $\pi_{n,t,\lambda}(K_n)$--invariant unit vector in
$\cH_{n,t,\lambda}$, then $\sigma(1K_{n,t}) = cu$ where $|c| = 1$ by
$K_{n,t}$--invariance, and then $\sigma = c\sigma_u$ by $K$--invariance.
In summary,
\begin{lemma}\label{k-invariant}
Let $t \in T$ and let $u$ be the unique {\rm (up to scalar multiple)}
$\pi_{n,t,\lambda}^\diamondsuit (K_{n,t})$--fixed unit vector in
$\cH_{n,t,\lambda}^\diamondsuit$  Then the section $\sigma_u$, given by
{\rm (\ref{induced-invariant})}, is the unique {\rm (up to scalar multiple)}
$\pi_{n,t,\lambda}(K)$--fixed unit vector in $\cH_{n,t,\lambda}$.
\end{lemma}

By Theorem \ref{iso-gelfand} we can apply (\ref{lim-diamond-induced})
and (\ref{lim-diamond-induced-1})
to the function spaces $\cE_{n,t,\lambda}^\diamond =
\cH_{n,t,\lambda}^\diamond \boxtimes (\cH_{n,t,\lambda}^\diamond)^*$
on the groups $G_{n,t} = N_n\rtimes K_{n,t}$.  Making use of
Lemma \ref{k-invariant} we have
\begin{proposition}\label{restrict-invariants}
If orthogonal projection
$\cE_{m,t,\lambda}^\diamondsuit \to \cE_{n,t,\lambda}^\diamondsuit$ sends
a nonzero right $K_{m,t}$--invariant function to a nonzero
right $K_{n,t}$--invariant function, then orthogonal projection
$\cE_{m,t,\lambda} \to \cE_{n,t,\lambda}$ sends a nonzero right
$K_m$--invariant function to a nonzero right $K_n$--invariant function.
\end{proposition}

E. Vinberg classified the maximal irreducible nilpotent Gelfand pairs.  
See \cite{V1}, \cite{V2}, or see \cite[Table 13.4.1]{W3}.  A Gelfand pair
$(G_n,K_n)$ is called {\sl maximal} if it is not obtained from another
Gelfand pair $(G'_n,K'_n)$ by the construction
$(G_n,K_n) = (G'_n/C,K'_n/(K'_n \cap C))$ for any nontrivial closed
connected central subgroup $C$ of $G'_n$.  And $(G_n,K_n)$ is called
{\sl irreducible} if $\Ad(K_n)$ is irreducible on $\gv_n = \gn_n / \gz$.
See \cite{W3} for the notation.  

{\small
\begin{equation} \label{vin-table}
\begin{tabular}{|r|c|c|c|c|c|}\hline
\multicolumn{6}{| c |}{Maximal Irreducible Nilpotent Gelfand
        Pairs $(N_n\rtimes K_n,K_n)$ \quad (\cite{V1}, \cite{V2})}\\
\hline \hline
 & Group $K_n$ & $\gv_n$ & $\gz$ &
   $\begin{smallmatrix} U(1) \text{ is} \\ \text{needed if}\end{smallmatrix}$ &
   $\begin{smallmatrix} \text{ max }\\\text{requires}\end{smallmatrix}$
   \\ \hline
1 & $SO(n)$ & $\R^n$ & $\Skew\R^{n\times n} = \gs\go(n)$ &  & \\ \hline
2 & $Spin(7)$ & $\R^8 = \O$ & $\R^7 = \Im\O$  &  & \\ \hline
3 & $G_2$ & $\R^7 = \Im\O$ & $\R^7 = \Im\O$ &  & \\ \hline
4 & $U(1)\cdot SO(n)$ & $\C^n$ & $\Im\C$ & & $n\ne 4$ \\ \hline
5 & $(U(1)\cdot) SU(n)$ & $\C^n$ & $\Lambda^2\C^n \oplus\Im\C$ &
        $n$ odd &  \\ \hline
6 & $SU(n), n$ odd & $\C^n$ & $\Lambda^2\C^n$ &  & \\ \hline
7 & $SU(n), n$ odd & $\C^n$ & $\Im\C$ &  & \\ \hline
8 & $U(n)$ & $\C^n$ & $\Im \C^{n\times n} = \gu(n)$ &  & \\ \hline
9 & $(U(1)\cdot) Sp(n)$ & $\H^n$ & $\Re \H^{n \times n}_0 \oplus \Im\H$ &  &
        \\ \hline
10 & $U(n)$ & $S^2\C^n$ & $\R$ & & \\ \hline
11 & $(U(1)\cdot) SU(n), n \geqq 3$ & ${\Lambda}^2\C^n$ & $\R$ & $n$ even &
        \\ \hline
12 & $U(1)\cdot Spin(7)$ & $\C^8$ & $\R^7 \oplus \R$ & & \\ \hline
13 & $U(1)\cdot Spin(9)$ & $\C^{16}$ & $\R$ & & \\ \hline
14 & $(U(1)\cdot) Spin(10)$ & $\C^{16}$ & $\R$ & & \\ \hline
15 & $U(1)\cdot G_2$ & $\C^7$ & $\R$ & & \\ \hline
\multicolumn{6}{r}{.... table continued on next page}
\end{tabular}
\end{equation}
}
{\normalsize
\addtocounter{equation}{-1}
\begin{equation}
\begin{tabular}{|r|c|c|c|c|c|}
\multicolumn{6}{ l }{.... table continued from previous page} \\ \hline
16 & $U(1)\cdot E_6$ & $\C^{27}$ & $\R$ & & \\ \hline
17 & $Sp(1)\times Sp(n)$ & $\H^n$ & $\Im \H = \gs\gp(1)$ & & $n \geqq 2$ \\ \hline
18 & $Sp(2)\times Sp(n)$ & $\H^{2\times n}$ &
        $\Im \H^{2\times 2} = \gs\gp(2)$ & & \\ \hline
19 & $(U(1)\cdot) SU(m) \times SU(n)$ &  &  &  & \\
   & $m,n \geqq 3$ & $\C^m\otimes \C^n$ & $\R$ & $m=n$ &   \\ \hline
20 & $(U(1)\cdot) SU(2) \times SU(n)$ & $\C^2 \otimes \C^n$ &
        $\Im \C^{2\times 2} = \gu(2)$ & $n=2$ & \\ \hline
21 & $(U(1)\cdot) Sp(2) \times SU(n)$ & $\H^2\otimes \C^n$ & $\R$ &
        $n \leqq 4$ & $n \geqq 3$ \\ \hline
22 & $U(2)\times Sp(n)$ & $\C^2 \otimes \H^n$ & $\Im \C^{2\times 2} = \gu(2)$ &
        & \\ \hline
23 & $U(3)\times Sp(n)$ & $\C^3 \otimes \H^n$ & $\R$ & & $n \geqq 2$
        \\ \hline
\end{tabular}
\end{equation}
}

Splitting some cases to retain parabolicity of $\{K_n\}$,
the strict direct systems in Table \ref{vin-table}, with $\dim \gz_n$
bounded, are

{\normalsize
\begin{equation} \label{ind-vin-table}
\begin{tabular}{|r|c|c|c|c|c|}\hline
\multicolumn{6}{| c |}{Direct Systems of Maximal Irreducible Nilpotent Gelfand
        Pairs $(N_n\rtimes K_n,K_n)$}\\
\hline \hline
 & Group $K_n$ & $\gv_n$ & $\gz_n$ &
   $\begin{smallmatrix} U(1) \text{ is} \\ \text{needed if}\end{smallmatrix}$ &    $\begin{smallmatrix} \text{ max }\\\text{requires}\end{smallmatrix}$
   \\ \hline
4a & $U(1)\cdot SO(2n)$ & $\C^{2n}$ & $\Im\C$ & & $n\ne 2$ \\ \hline
4b & $U(1)\cdot SO(2n+1)$ & $\C^{2n+1}$ & $\Im\C$ & &  \\ \hline
7 & $SU(n), n$ odd & $\C^n$ & $\Im\C$ &  & \\ \hline
10 & $U(n)$ & $S^2\C^n$ & $\R$ & & \\ \hline
11 & $(U(1)\cdot) SU(n), n \geqq 3$ & ${\Lambda}^2\C^n$ & $\R$ & $n$ even &
        \\ \hline
17 & $Sp(1)\times Sp(n)$ & $\H^n$ & $\Im \H = \gs\gp(1)$ & & $n \geqq 2$
        \\ \hline
18 & $Sp(2)\times Sp(n)$ & $\H^{2\times n}$ &
        $\Im \H^{2\times 2} = \gs\gp(2)$ & & \\ \hline
19 & $(U(1)\cdot) SU(m) \times SU(n)$ &  &  &  & \\
   & $m,n \geqq 3$ & $\C^m\otimes \C^n$ & $\R$ & $m=n$ &   \\ \hline
20a & $SU(2) \times SU(n), n \geqq 3$ & $\C^2 \otimes \C^n$ &
        $\Im \C^{2\times 2} = \gu(2)$ &  & \\ \hline
20b & $U(2) \times SU(n)$ & $\C^2 \otimes \C^n$ &
        $\Im \C^{2\times 2} = \gu(2)$ &  & \\ \hline
21 & $(U(1)\cdot) Sp(2) \times SU(n)$ & $\H^2\otimes \C^n$ & $\R$ &
        $n \leqq 4$ & $n \geqq 3$ \\ \hline
22 & $U(2)\times Sp(n)$ & $\C^2 \otimes \H^n$ & $\Im \C^{2\times 2} = \gu(2)$ &
        & \\ \hline
23 & $U(3)\times Sp(n)$ & $\C^3 \otimes \H^n$ & $\R$ & & $n \geqq 2$
        \\ \hline
\end{tabular}
\end{equation}}

In each case of Table \ref{ind-vin-table}, \cite[Theorem 14.4.3]{W3} says that
$N_n$ has square integrable representations.
In the cases $\dim \gz > 1$ of Table \ref{ind-vin-table} we have
$K_n = K^\dagger_n\cdot L$ where the big factor $K_n^\dagger$ acts trivially
on $\gz$ and the small factor $L$ acts on $\gz$ by its adjoint
representation.  Summarizing these observations,

\begin{proposition} Each direct system $\{(G_n,K_n)\}$ of
{\rm Table \ref{ind-vin-table}} has the properties {\rm (i)} $\{K_n\}$ is
parabolic, {\rm (ii)} the $\{K_{n,t}\}$ are parabolic for each $t \in T$,
and {\rm (iii)} $N_n$ has square integrable representations.
\end{proposition}

From Theorem \ref{iso-gelfand} and the argument of Lemma \ref{kntriv},
\begin{corollary}\label{kntriv-vin}
Let $\{(G_n,K_n)\}$ be a direct system of
{\rm Table \ref{ind-vin-table}} and let $t \in T$.
Then $\pi_{n,t,\lambda}^\diamondsuit$ has a
nonzero $K_{n,t}$--fixed vector if and only if $\kappa_{n,t,\lambda}^*$
occurs as a subrepresentation of $\widetilde{\pi_{n,t}}|_{K_{n,t}}$, and in
that case the space of $K_{n,t}$--fixed vectors has dimension $1$.
\end{corollary}

Now everything goes essentially as in the Heisenberg nilmanifold cases
of Section \ref{sec4}.  
We have isometric $G_n$--equivariant injections
\begin{equation}\label{k-zeta-nilp7}
\widetilde{\zeta}_{m,n,t,\lambda}: \cE_{n,t,\lambda}^{K_n} 
\to \cE_{m,t,\lambda}^{K_m} \text{ by }
\widetilde{\zeta}_{m,n,t,\lambda}((\deg\pi_{n,t,\lambda})^{1/2}f)
=c_{m,n,t,\lambda}(\deg\pi_{m,t,\lambda})^{1/2} f
\end{equation}
where as in (\ref{alpha1}), $c_{m,n,t,\lambda}$ is the ratio (\ref{def-c})
of lengths of $K_n$, $K_m$ fixed unit vectors.
Integrating on $t$ and summing $\lambda$ gives
isometric $G_n$--equivariant injections
$\widetilde{\zeta}_{m,n}: L^2(G_n/K_n) \to L^2(G_m/K_m)$, as follows.

\begin{theorem}\label{lim-sd7}
For each of the direct systems of {\rm Table \ref{ind-vin-table}}
denote $G = \varinjlim \{G_n\}$, $N = \varinjlim N_n$ and
$K = \varinjlim \{K_n\}$.  Note $G = N \rtimes K$.
Then $\{L^2(G_n/K_n),\widetilde{\zeta}_{m,n}\}$ is a strict direct system 
of Hilbert spaces
in which the  maps $\widetilde{\zeta}_{m,n}: L^2(G_n) \to L^2(G_m)$ are
$G_n$--equivariant injections.
That gives us the left regular representation of $G$
on the Hilbert space $L^2(G): = \varinjlim \{L^2(G_n), \zeta_{m,n}\}$.
Further, that left/right regular representation
is the multiplicity--free
$\sum_\lambda \int_{\Ad^*(L)\backslash \gz^*} (\pi_{\Ad^*(k)t,\Ad^*(k)\lambda}
\boxtimes \pi_{\Ad^*(k)t,\Ad^*(k)\lambda}^*)\,\,dk$ where
$\pi_{t,\lambda} := \varinjlim \pi_{n,t,\lambda}$.
\end{theorem}

Similarly using 
$w_{m,\lambda} = c_{m,n,t,\lambda}w_{n,\lambda} + x$ with
$x \perp \cE_{n,t,\lambda}$ with $0 < c_{n,t,\lambda} \leqq 1$,
\begin{equation}\label{reg-fin-heis-kk}
 \cA(G_n/K_n) := \cA(G_n)^{K_n} =
        \cA(G_n)\cap \cE_{n,t,\lambda}^{K_n}
\end{equation}
leading to direct systems and their limits by assembling the maps
\begin{equation}\label{reg-fin-nil-k}
\nu_{m,n,t,\lambda}:\cA(\pi_{n,t,\lambda})^{K_n} \hookrightarrow
  \cA(\pi_{m,t,\lambda})^{K_m} \text{ by }
        f_{u,v_n,n} \mapsto f_{u,v_m,m} \text{ where }
        p_{m,n,t,\lambda}(v_m) = v_n
\end{equation}
Then we have direct systems and their limits
\begin{equation}\label{reg-inf-nil-k}
\begin{aligned}
& \cA(G) = \varinjlim \{\cA(G_n), \zeta_{m,n}\}
\text{ where } \zeta_{m,n,t,\lambda} : \cE_{n,t,\lambda}
\hookrightarrow \cE_{m,t,\lambda} \text{  (\ref{lim-diamond-induced}), and}\\
&\cA(G/K) =
\varinjlim \{\cA(G_n)/K_n), \nu_{m,n}\}
\text{ where } \nu_{m,n,t,\lambda} : \cE_{n,t,\lambda}^{K_n}
\hookrightarrow \cE_{m,t,\lambda}^{K_m} \text{  (\ref{reg-fin-nil-k})}.
\end{aligned}
\end{equation}

As before, each $\cA(G_n/K_n)$ is dense in $L^2(G_n/K_n)$, and we
pass this comparison to the limit with the maps
\begin{equation}\label{kkkk-heis-rel-map-sys}
\widetilde{\eta}_{n,t,\lambda}: \cA(\pi_{n,t,\lambda})^{K_n} \to
\cH_{\pi_{n,t,\lambda}}\otimes(w_{n,t,\lambda^*}\C) \text{ by }
f \mapsto c_{n,1,t,\lambda}\,\,\sqrt{\deg \pi_{n,t,\lambda}}\, f.
\end{equation}
We conclude
\begin{proposition}\label{kkkk-heis-quo-comparison}
The $\widetilde{\eta}_{n,t,\lambda}$
satisfy $(\widetilde{\eta}_{m,t,\lambda}\circ
\nu_{m,n,t,\lambda})(f) =
(\widetilde{\zeta}_{m,n,\lambda}\circ\widetilde{\eta}_{n,t,\lambda})(f)$
for $f \in \cA(\pi_{n,t,\lambda})^{K_n}$  Thus they inject
the direct system $\{\cA(G_n/K_n), \nu_{m,n}\}$
of regular functions
into the direct system $\{L^2(G_n/K_n),\widetilde{\zeta}_{m,n}\}$.
That map of direct systems defines a $G$--equivariant
injection
$
\widetilde{\eta}: \cA(G/K) \to L^2(G/K)
$
with dense image.  In particular $\widetilde{\eta}$ defines a pre Hilbert space
structure on $\cA(G/K)$ with completion isometric
to $L^2(G/K)$.
\end{proposition}

\subsection{Reducible Quotients.}\label{sec5d}
There are more many strict direct sequences of nilmanifold Gelfand pairs
$(G_n,K_n)$, e.g. those for which the action of $K_n$ on $\gn_n/\gz_n$ is
reducible.  These $(G_n/K_n)$ are constructed from certain basic ones that
satisfy several technical conditions (indecomposable, principal, maximal
and $Sp(1)$--saturated).  See \cite{Y1},\cite{Y2}, \cite{Y3}, \cite{W3} and 
\cite{W5}.  The basic such direct systems, with $K_n$ reducible on
$\gn_n/\gz_n$, $\dim\gz_n$ bounded and $\{K_n\}$ parabolic, are tabulated 
in \cite[Table 9.15]{W5} as follows.  Here the numbering comes from
\cite[Table 9.14]{W5}, $N_\ell = N_\ell'\times Z_\ell$ with $Z_\ell$
central and maximal for that, $G_\ell = N_\ell\rtimes K_\ell$ and
$G'_\ell = N'_\ell\rtimes K_\ell$\,.  The generalized Heisenberg algebra 
$\gh_{n,\F}$ is $\Im\F + \F^n$ with $[(z,w),(z',w')] = (z+z'+\Im
\langle w,w'\rangle, w+w')$ where $\F = \C$, $\H$ (quaternions) or
$\O$ (octonions).  

{\small
\begin{equation} \label{ind-vin-table-ipms}
\begin{tabular}{|c|l|l|l|l|l|}\hline
\multicolumn{6}{| c |}{Strict Direct Systems $\{(G_\ell,K_\ell)\}$
and $\{(G'_\ell,K'_\ell)\}$
of Gelfand Pairs with
$\dim \gz'_\ell$ Bounded}\\
\hline \hline
 & Group $K_\ell$ & $K_\ell$--module $\gv_\ell$ 
	& $\gz_\ell' = [\gn_\ell,\gn_\ell]$ &
        $\gz''_\ell$ & Algebra $\gn'_\ell$
        \\ \hline
1 & $U(n)$ & $\C^n$ & $\R$ & $\gs\gu(n)$ &
        $\gh_{n;\C}$  \\ \hline
3 & $U(1)\times U(n)$ & $\C^n \oplus \Lambda^2\C^n$ & $\R \oplus \R$ & $0$ &
        $\gh_{n;\C} \oplus \gh_{n(n-1)/2;\C}$ \\ \hline
6 & $S(U(4)\times U(m))$ & $\C^{4 \times m}$ & $\R$ & $\R^6$ &
        $\gh_{4m;\C}$ \\ \hline
7 & $U(m)\times U(n)$ & $\C^{m\times n}\oplus \C^m$ & $\R \oplus \R$ & $0$ &
        $\gh_{mn;\C} \oplus \gh_{m;\C}$ \\ \hline
8 & $U(1)\times Sp(n)\times U(1)$ & $\C^{2n}\oplus \C^{2n}$ & $\R\oplus \R$ &
         $0$ & $\gh_{2n;\C} \oplus \gh_{2n;\C}$ \\ \hline
9 & $Sp(1)\times Sp(n)\times U(1)$ & $\H^n \oplus \H^n$ & $\Im \H \oplus \R$ &
         $0$ & $\gh_{n;\H} \oplus \gh_{2n;\C}$ \\ \hline
10 & $Sp(1)\times Sp(n)\times Sp(1)$ & $\H^n\oplus\H^n$ & $\Im\H\oplus\Im\H$ &
         $0$ & $\gh_{n;\H} \oplus \gh_{n;\H}$ \\ \hline
11a & $Sp(n)\times Sp(1)\times Sp(m)$ & $\H^n$ & $\Im\H$ &
        $\H^{n\times m}$ & $\gh_{n;\H}$ \\ \hline
11b & $Sp(n)\times U(1)\times Sp(m)$ & $\H^n$ & $\Im\H$ &
        $\H^{n\times m}$ & $\gh_{n;\H}$ \\ \hline
11c & $Sp(n)\times \{1\} \times Sp(m)$ & $\H^n$ & $\Im\H$ &
        $\H^{n\times m}$ & $\gh_{n;\H}$ \\ \hline
18a & $SU(n) \times SU(2)$ &
        $\C^{n\times 2}$ & $\R$ & $\gs\gu(2)$ &
        $\gh_{2n;\C}$ \\ \hline
18b & $U(n) \times SU(2)$ &
        $\C^{n\times 2}$ & $\R$ & $\gs\gu(2)$ &
        $\gh_{2n;\C}$ \\ \hline
18c & $U(1)Sp(\tfrac{n}{2}) \times SU(2)$ &
        $\C^{n\times 2}$ & $\R$ & $\gs\gu(2)$ &
        $\gh_{2n;\C}$ \\ \hline
19a & $SU(n) \times U(2)$ &
        $\C^{n\times 2}\oplus\C^2$ & $\R\oplus\R$ & $0$ &
        $\gh_{2n;\C} \oplus \gh_{2;\C}$ \\ \hline
19b & $U(n) \times U(2)$ &
        $\C^{n\times 2}\oplus\C^2$ & $\R\oplus\R$ & $0$ &
        $\gh_{2n;\C} \oplus \gh_{2;\C}$ \\ \hline
19c & $U(1)Sp(\tfrac{n}{2}) \times U(2)$ &
        $\C^{n\times 2}\oplus\C^2$ & $\R\oplus\R$ & $0$ &
        $\gh_{2n;\C} \oplus \gh_{2;\C}$ \\ \hline
20aa & $SU(n) \times SU(2)\times SU(m)$ &
        $\C^{n\times 2}\oplus\C^{2\times m}$ &
        $\R \oplus \R$ & $0$ & $\gh_{2n;\C} \oplus \gh_{2m;\C}$ \\ \hline
20ab & $SU(n) \times SU(2)\times U(m)$ &
        $\C^{n\times 2}\oplus\C^{2\times m}$ &
        $\R \oplus \R$ & $0$ & $\gh_{2n;\C} \oplus \gh_{2m;\C}$ \\ \hline
20ac & $SU(n) \times SU(2)\times U(1)Sp(\tfrac{m}{2})$ &
        $\C^{n\times 2}\oplus\C^{2\times m}$ &
        $\R \oplus \R$ & $0$ & $\gh_{2n;\C} \oplus \gh_{2m;\C}$ \\ \hline
20ba & $U(n) \times SU(2)\times SU(m)$ &
        $\C^{n\times 2}\oplus\C^{2\times m}$ &
        $\R \oplus \R$ & $0$ & $\gh_{2n;\C} \oplus \gh_{2m;\C}$ \\ \hline
20bb & $U(n) \times SU(2)\times U(m)$ &
        $\C^{n\times 2}\oplus\C^{2\times m}$ &
        $\R \oplus \R$ & $0$ & $\gh_{2n;\C} \oplus \gh_{2m;\C}$ \\ \hline
20bc & $U(n) \times SU(2)\times U(1)Sp(\tfrac{m}{2})$ &
        $\C^{n\times 2}\oplus\C^{2\times m}$ &
        $\R \oplus \R$ & $0$ & $\gh_{2n;\C} \oplus \gh_{2m;\C}$ \\ \hline
20ca & $U(1)Sp(\tfrac{n}{2}) \times SU(2)\times SU(m)$ &
        $\C^{n\times 2}\oplus\C^{2\times m}$ &
        $\R \oplus \R$ & $0$ & $\gh_{2n;\C} \oplus \gh_{2m;\C}$ \\ \hline
20cb & $U(1)Sp(\tfrac{n}{2}) \times SU(2)\times U(m)$ &
        $\C^{n\times 2}\oplus\C^{2\times m}$ &
        $\R \oplus \R$ & $0$ & $\gh_{2n;\C} \oplus \gh_{2m;\C}$ \\ \hline
20cc & $U(1)Sp(\tfrac{n}{2}) \times SU(2)\times$ &
        $\C^{n\times 2}\oplus\C^{2\times m}$ &
        $\R \oplus \R$ & $0$ & $\gh_{2n;\C} \oplus \gh_{2m;\C}$ \\
   & \phantom{$\{SU(m),U(m)\}$} $U(1)Sp(\tfrac{m}{2})\}$ & & & & \\ \hline
21a & $SU(n) \times SU(2) \times U(4)$ &
        $\C^{n\times 2}\oplus\C^{2\times 4}$ &
        $\R \oplus \R$ & $\R^6$ &
        $\gh_{2n;\C} \oplus \gh_{8;\C}$ \\ \hline
21b & $U(n) \times SU(2) \times U(4)$ &
        $\C^{n\times 2}\oplus\C^{2\times 4}$ &
        $\R \oplus \R$ & $\R^6$ &
        $\gh_{2n;\C} \oplus \gh_{8;\C}$ \\ \hline
21c & $U(1)Sp(\tfrac{n}{2}) \times SU(2) \times U(4)$ &
        $\C^{n\times 2}\oplus\C^{2\times 4}$ &
        $\R \oplus \R$ & $\R^6$ &
        $\gh_{2n;\C} \oplus \gh_{8;\C}$ \\ \hline
\end{tabular}
\end{equation}
}
One obtains the structure of $L^2(G/K)$ and $\cA(G/K)$ 
for the cases of Table \ref{ind-vin-table-ipms}
by a straightforward modification of the considerations involved for the
direct systems of Table \ref{ind-vin-table}.  We leave the details to the
reader.

\vskip .2 in

\noindent
\phantom{XXXX}Department of Mathematics\hfill\newline\noindent
\phantom{XXXX}University of California\hfill\newline\noindent
\phantom{XXXX}Berkeley, CA 94720--3840, USA \hfill\newline\noindent
\hfill\newline\noindent
\phantom{XXXX}{\tt jawolf@math.berkeley.edu}

\end{document}